\documentclass[reqno]{amsart}
\usepackage{graphicx}
\usepackage{cite}
\usepackage{amssymb}
\usepackage{tikz}
\newtheorem{thm}{Theorem}[section]
\newtheorem{cor}[thm]{Corollary}
\newtheorem{lem}[thm]{Lemma}

\theoremstyle{definition}
\newtheorem{defn}[thm]{Definition}
\theoremstyle{remark}
\newtheorem{rem}[thm]{Remark}
\numberwithin{equation}{section}

\begin{document}

\title[Orthogonal trigonometric polynomials from Riemann-Hilbert view]{Orthogonal trigonometric polynomials from Riemann-Hilbert view}%
\author{Zhihua Du}%


\begin{abstract}
In this work, some theorems are established for orthogonal trigonometric polynomials (OTP) including Favard, Baxter, Geronimus,
Rakhmanov, Szeg\"o and the strong Szeg\"o theorems which are important in the
theory of orthogonal polynomials on the unit circle (OPUC). All these results are based on the mutual representation theorem of OPUC and OTP which deduced by a Riemann-Hilbert problem simultaneously characterizing them. In addition, Szeg\"o recursions, four-term recurrences and some new identities  for OPUC and OTP are also obtained by using their Riemann-Hilbert characterizations respectively.
\end{abstract}
\maketitle
\section{Introduction}

It is well known that orthogonal polynomials have a great history and continuing important applications in mathematics, physics and engineering and so on \cite{sze,sim1,sim2,deift,mehta}. Boundary value problems for analytic functions is a living research field with a beautiful and rich theory as well as diverse and interesting applications. It also has a fascinating history which can be traced back to the origins of function theory and comes into sight after 1851 via Riemann's famous habilitation thesis. Riemann's treatment on these problems is heuristic. It was Hilbert who first proposed a partly rigorous approach to attack the problems in the linear case. The main defect of Hilbert's approach lies in the ignorance of the indexes of the problems which pointed out by F. Noether. For this reason, nowadays these problems are usually called Riemann-Hilbert problems (for short, RHPs). In the past almost thirty years, a remarkable fact is that one can construct some (usually, $2\times2$) matrix-valued RHPs to characterize many different types of orthogonal polynomials with respect to general weight functions or probability measures. These RHPs are always called Riemann-Hilbert characterizations (simply, RH characterizations; or more simply, RHCs) for the corresponding orthogonal polynomials.

In fact, RHPs appear in many different settings. There are many systematic approaches to formulate RH characterizations for some interesting problems in modern studies. Nevertheless, RH characterizations for orthogonal polynomials come ``out of the blue" according to Deift's view \cite{deift1}. In this regard, the first breakthrough was due to Fokas, Its and Kitaev \cite{fik}. There they proposed the RH characterization for orthogonal polynomials on the real line (simply, OPRL). More precisely, they formulated the following $2\times2$ matrix-valued RHP for a $2\times2$ matrix-valued function $\mathcal{Y}: \mathbb{C}\setminus \mathbb{R}\rightarrow \mathbb{C}^{2\times 2}$ satisfying
\begin{equation} (\mbox{RHC for OPRL})\,\,
\begin{cases}
\mathcal{Y}\,\, \mbox{is analytic in}\,\,\mathbb{C}\setminus \mathbb{R},\vspace{2mm}\\
\mathcal{Y}^{+}(x)=\mathcal{Y}^{-}(x)\left(
                   \begin{array}{cc}
                     1 & w(x) \\
                     0 & 1 \\
                   \end{array}
                 \right)
                 \,\,\mbox{for} \,\,x\in \mathbb{R},\vspace{2mm}\\
\mathcal{Y}(z)=\left(I+O\left(\frac{1}{z}\right)\right)\left(
                                    \begin{array}{cc}
                                      z^{n} & 0 \\
                                      0 & z^{-n} \\
                                    \end{array}
                                  \right)
                                  \,\,\mbox{as}\,\, z\rightarrow \infty,
\end{cases}
\end{equation}
where $w$ is a weight function on $\mathbb{R}$ and $I$ is the $2\times2$ identity matrix.

In \cite{bdj}, Baik, Deift and Johansson proposed the RH characterization for orthogonal polynomials on the unit circle (concisely, OPUC). That is, for a $2\times2$ matrix-valued function $Y: \mathbb{C}\setminus \partial \mathbb{D}\rightarrow \mathbb{C}^{2\times 2}$, the following conditions are fulfilled:
\begin{equation} (\mbox{RHC for OPUC})\,\,
\begin{cases}
Y\,\, \mbox{is analytic in}\,\,\mathbb{C}\setminus \partial \mathbb{D},\vspace{2mm}\\
Y^{+}(t)=Y^{-}(t)\left(
                   \begin{array}{cc}
                     1 & t^{-n}w(t) \\
                     0 & 1 \\
                   \end{array}
                 \right)
                 \,\,\mbox{for} \,\,t\in \partial \mathbb{D},\vspace{2mm}\\
Y(z)=\left(I+O(\frac{1}{z})\right)\left(
                                    \begin{array}{cc}
                                      z^{n} & 0 \\
                                      0 & z^{-n} \\
                                    \end{array}
                                  \right)
                                  \,\,\mbox{as}\,\, z\rightarrow \infty,
\end{cases}
\end{equation}
where $\mathbb{D}$ is the unit disc, $\partial \mathbb{D}$ is the unit circle, $w$ is a weight function on $\partial \mathbb{D}$ and $I$ is the $2\times2$ identity  matrix.

With respect to orthogonal trigonometric polynomials (simply, OTP), in \cite{dd08}, Du and the author constructed the RH characterization for them.  More precisely, it is the following $2\times2$ matrix-valued RHP: for a $2\times2$ matrix-valued function $\mathfrak{Y}: \mathbb{C}\setminus \partial \mathbb{D}\rightarrow \mathbb{C}^{2\times 2}$ satisfying
\begin{equation} (\mbox{RHC for OTP})\,\,
\begin{cases}
\mathfrak{Y}\,\, \mbox{is analytic in}\,\,\mathbb{C}\setminus \partial \mathbb{D},\vspace{2mm}\\
\mathfrak{Y}^{+}(t)=\mathfrak{Y}^{-}(t)\left(
                   \begin{array}{cc}
                     1 & t^{-2n}w(t) \\
                     0 & 1 \\
                   \end{array}
                 \right)
                 \,\,\mbox{for} \,\,t\in \partial \mathbb{D},\vspace{2mm}\\
\mathfrak{Y}(z)=\left(I+O(\frac{1}{z})\right)\left(
                                    \begin{array}{cc}
                                      z^{2n} & 0 \\
                                      0 & z^{-2n+1} \\
                                    \end{array}
                                  \right)
                                  \,\,\mbox{as}\,\, z\rightarrow \infty,\vspace{2mm}\\
\mathfrak{Y}_{11}(0)=\mathfrak{Y}_{21}(0)=0,
\end{cases}
\end{equation}
where $\mathbb{D}$ is the unit disc, $\partial \mathbb{D}$ is the unit circle, $w$ is a weight function on $\partial \mathbb{D}$ and $I$ is the $2\times2$ identity  matrix.

Observing all of the above RH characterizations, the innovation from OPRL to OPUC is that the 12 entry in the jump matrix with $t^{-n}w$ replacing $w$. However, from OPUC to OTP, there is another completely new innovation except the 12 entry $t^{-2n}w$ in place of $t^{-n}w$ in the jump matrix. That is, in the matrix about the growth conditions at $\infty$ for $Y$, the 11 entry $z^{n}$ is replaced by $z^{2n}$ and the 22 entry $z^{-n}$ is replaced by $z^{-2n+1}$. Moreover, the 11 and 21 entries are prescribed to be $0$ at the origin. Based on such innovations, a surprisingly remarkable fact is that the RHP (1.3) can be characterized for both OTP and OPUC. For this reason, Du and the author discovered and established the mutual representation theorem for OTP and OPUC. It becomes a bridge connecting such two isolated classes of orthogonal polynomials. However, it should be pointed out that the RHC (1.2) for OPUC can also be as a RHC for OTP when replacing $n$ with $2n-1$ (see Remark 3.2 in \cite{dd08}). Such two RHCs can be transformed to each other by an explicit $2\times2$ matrix-valued multiplier (see the uniqueness part of the proof of Theorem 3.1 in \cite{dd08}). Nevertheless, we still use the RHP (1.3) as the RH characterization for OTP in the present paper.

In addition, for general orthogonal polynomials (simply, GOP), the author also formulated a semi-conjugate $2\times2$ matrix-valued boundary value problem to characterize orthogonal polynomials on an arbitrary smooth Jordan curve in $\mathbb{C}$ (see \cite{d}). But it isn't a RHP since the semi-conjugate operator appears. More precisely, for a $2\times2$ matrix-valued function $\mathrm{Y}: \mathbb{C}\setminus \Gamma\rightarrow \mathbb{C}^{2\times 2}$, the following conditions are satisfied:
\begin{equation} (\mbox{RHC for GOP})\,\,
\begin{cases}
\mathrm{Y}\,\, \mbox{is analytic in}\,\,\mathbb{C}\setminus \Gamma,\vspace{2mm}\\
(\mathrm{D}\mathrm{Y})^{+}(t)=(\mathrm{D}\mathrm{Y})^{-}(t)\left(
                   \begin{array}{cc}
                     1 & w(t)s^{\prime}(t) \\
                     0 & 1 \\
                   \end{array}
                 \right)
                 \,\,\mbox{for} \,\,t\in \Gamma,\vspace{2mm}\\
\mathrm{Y}(z)=\left(I+O(\frac{1}{z})\right)\left(
                                    \begin{array}{cc}
                                      z^{n} & 0 \\
                                      0 & z^{-n} \\
                                    \end{array}
                                  \right)
                                  \,\,\mbox{as}\,\, z\rightarrow \infty,
\end{cases}
\end{equation}
where $\Gamma$ is an arbitrary smooth Jordan curve in $\mathbb{C}$ oriented counter-clockwisely, $s(t)$ is the arc-length function, $w$ is a weight function on $\Gamma$, $I$ is the $2\times2$ identity  matrix and the semi-conjugate operator $\mathrm{D}$ is defined by
\begin{equation}
    ({\mathrm{D}}\mathrm{Y})(z)=\left(
              \begin{array}{ll}
              \overline{\mathrm{Y}_{11}(z)} \hspace{2mm} \mathrm{Y}_{12}(z) \\
              \overline{\mathrm{Y}_{21}(z)} \hspace{2mm} \mathrm{Y}_{22}(z)
              \end{array}
           \right)
            \quad {\mathrm{for}} \quad \mathrm{Y}(z)=
            \left(
            \begin{array}{ll}
              \mathrm{Y}_{11}(z) \hspace{2mm} \mathrm{Y}_{12}(z) \\
              \mathrm{Y}_{21}(z) \hspace{2mm} \mathrm{Y}_{22}(z)
            \end{array}
            \right),
\end{equation}
in which $\overline{z}$ is the conjugate complex number of $z$.

As some applications of the mutual representation theorem, four-term recurrences, Christoffel-Darboux formulae and some properties of zeros for OTP were obtained in \cite{dd08}. In fact, by the mutual representation theorem, some important theorems (such as Favard, Baxter, Geronimus, Rakhmanov, Szeg\"o and the
strong Szeg\"o) in the theory of OPUC can be established for OTP. This is one of the themes in the present paper.

At present, together with the nonlinear steepest descent method due to Deift and Zhou \cite{dz}, Riemann-Hilbert problems are mainly applied into some asymptotic analysis problems on integrable systems, orthogonal polynomials, combinatorics and random matrices, etc. \cite{bdj,deift,mehta}. However, except for asymptotic analysis \cite{bdj,dd06,dd08,dz}, Riemann-Hilbert problems can also be used to some analytic and algebraic problems such as yielding some difference equations, differential equations and so on \cite{deift,deift1}. As an example to show the subtlety and power of RHPs in this facet,  a new proof is given for Szeg\"o recursions of OPUC and four-term recurrences of OTP by using their RH characterizations in Section 4. As a byproduct, some new identities on Cauchy integrals for both OPUC and OTP as well as Hilbert transforms for OTP are also obtained.

This paper is organized as follows.

In Section 2, some definitions and notations are introduced involving OPUC and OTP as well some coefficients about them such as Verblunsky coefficients and so on, the mutual representation theorem and some of its consequences are given.

In Section 3, some theorems for OTP are obtained in terms of the mutual representation theorem including Favard, Baxter, Geronimus,
Rakhmanov, Szeg\"o and the strong Szeg\"o theorems which are important in the theory of OPUC. However, Favard theorem in this section is only in a weak form.

As stated above, in Section 4, some identities such as Szeg\"o recursions of OPUC, four-term recurrences of OTP and so on are obtained by using RH characterizations for OPUC and OTP respectively.

The final section is mainly devoted to prove a Favard theorem stronger than the one in Section 3.

\section{Mutual representation and its consequences}
Let $\mathbb{D}$ be the unit disc in the complex plane,
$\partial\mathbb{D}$ be the unit circle and $\mu$ be a nontrivial
probability measure on $\partial \mathbb{D}$ (i.e. with infinity
support, nonnegative and $\mu(\partial \mathbb{D})=1$). Throughout this paper,
by decomposition, we always write
\begin{equation}
d\mu(\tau)=w(\tau)\frac{d\tau}{2\pi i\tau}+d\mu_{s}(\tau),
\end{equation}
where $\tau\in\partial \mathbb{D}$, $w(\tau)=2\pi i\tau
d\mu_{ac}/d\tau$ in which $d\mu_{ac}$ is the absolutely continuous part of $d\mu$, and $d\mu_{s}$ is the singular part of $d\mu$.

Introduce two class of inner products, one is complex as follows
\begin{equation}
\langle f, g\rangle_{\mathbb{C}}=\int_{\partial
\mathbb{D}}\overline{f(\tau)}g(\tau)d\mu(\tau)
\end{equation}
with norm $||f||_{\mathbb{C}}=[\int_{\partial
\mathbb{D}}|f(\tau)|^{2}d\mu(\tau)]^{1/2}$, where $f, g$ are complex
integrable functions on $\partial \mathbb{D}$. The other is real and defined by
\begin{equation}
\langle f, g\rangle_{\mathbb{R}}=\int_{\partial
\mathbb{D}}f(\tau)g(\tau)d\mu(\tau)
\end{equation}
with norm $||f||_{\mathbb{R}}=[\int_{\partial
\mathbb{D}}|f(\tau)|^{2}d\mu(\tau)]^{1/2}$, where $f, g$ are real
integrable functions on $\partial \mathbb{D}$.

By the complex inner product (2.2), applying Gram-Schmidt procedure
to the following system
\begin{equation*}
\{1,z,z^{2},\ldots,z^{n},\ldots\},
\end{equation*}
where $z\in \mathbb{C}$, we get the unique system $\{\Phi_{n}(z)\}$
of monic orthogonal polynomials on the unit circle with respect to
$\mu$ satisfying
\begin{equation}
\langle\Phi_{n},\Phi_{m}\rangle_{\mathbb{C}}=\kappa_{n}^{-2}\delta_{nm}\,\,\,\text{with}\,\,\,
\kappa_{n}>0.
\end{equation}
Then the orthonormal polynomials $\varphi_{n}(z)$ on the unit circle
satisfy
\begin{equation}
\langle\varphi_{n},\varphi_{m}\rangle_{\mathbb{C}}=\delta_{nm}\,\,\,\text{and}\,\,\,
\varphi_{n}(z)=\kappa_{n}\Phi_{n}(z).
\end{equation}

For any polynomial $Q_{n}$ of order $n$, its reversed polynomial
$Q_{n}^{*}$ is defined by
\begin{equation}
Q_{n}^{*}(z)=z^{n}\overline{Q_{n}(1/\overline{z})}.
\end{equation}

One famous property of OPUC is Szeg\"o recurrence \cite{sze}, i.e.
\begin{equation}
\Phi_{n+1}(z)=z\Phi_{n}(z)-\overline{\alpha}_{n}\Phi^{*}_{n}(z),
\end{equation}
where $\alpha_{n}=-\overline{\Phi_{n+1}(0)}$ are called Verblunsky
coefficients. It is well known that $\alpha_{n}\in \mathbb{D}$ for
$n\in \mathbb{N}\cup\{0\}$. By convention, $\alpha_{-1}=-1$ (see
\cite{sim1}). Szeg\"o recurrence (2.7) is extremely useful in the
theory of OPUC. Especially, Verblunsky coefficients play an
important role in many interesting problems for OPUC (see \cite{sim1,sim2}).

Using the real inner product (2.3) and Gram-Schmidt procedure to the
following over $\mathbb{R}$ linearly independent ordered set
  \begin{equation}
\Big\{1, \frac {z-z^{-1}}{2i}, \frac {z+z^{-1}}{2},  \ldots,
  \frac{z^{n}-z^{-n}}{2i}, \frac{z^{n}+z^{-n}}{2},  \ldots  \Big\},
  \end{equation}
where $z\in \mathbb{C}\setminus\{0\}$, we get the unique system
\begin{equation}
\{1, b_{1}\pi_{1}(z), a_{1}\sigma_{1}(z), \ldots, b_{n}\pi_{n}(z),
a_{n}\sigma_{n}(z), \ldots\}
\end{equation} of
the ``monic" orthogonal Laurent polynomials (concisely, OLP) of the first class on the unit circle
with respect to $\mu$ fulfilling
\begin{equation}
\langle\pi_{m},\sigma_{n}\rangle_{\mathbb{R}}=0,
\langle\pi_{m},\pi_{n}\rangle_{\mathbb{R}}=\langle\sigma_{m},\sigma_{n}\rangle_{\mathbb{R}}=\delta_{mn},\,\,\,m,n=1,2,\ldots
\end{equation}
and
\begin{equation}
a_{n}\sigma_{n}(z)=\frac{z^{n}+z^{-n}}{2}-\beta_{n}b_{n}\pi_{n}(z)-\imath_{n}a_{n-1}\sigma_{n-1}(z)
-\jmath_{n}b_{n-1}\pi_{n-1}(z)+\text{lower order}
\end{equation}
as well as
\begin{equation}
b_{n}\pi_{n}(z)=\frac{z^{n}-z^{-n}}{2i}-\varsigma_{n}a_{n-1}\sigma_{n-1}(z)
-\zeta_{n}b_{n-1}\pi_{n-1}(z)+\text{lower order},
\end{equation}
where $a_{n},b_{n}>0$, which are respectively the norms of the
``monic" orthogonal Laurent polynomials of the first class given by right hand sides
of (2.11) and (2.12),
\begin{equation}
\beta_{n}=\langle\frac{z^{n}+z^{-n}}{2},b_{n}^{-1}\pi_{n}\rangle_{\mathbb{R}},
\end{equation}
\begin{equation}
\imath_{n}=\langle\frac{z^{n}+z^{-n}}{2},a_{n-1}^{-1}\sigma_{n-1}\rangle_{\mathbb{R}},\,\,
\jmath_{n}=\langle\frac{z^{n}+z^{-n}}{2},b_{n-1}^{-1}\pi_{n-1}\rangle_{\mathbb{R}}
\end{equation}
and
\begin{equation}
\varsigma_{n}=\langle\frac{z^{n}-z^{-n}}{2i},a_{n-1}^{-1}\sigma_{n-1}\rangle_{\mathbb{R}},\,\,
\zeta_{n}=\langle\frac{z^{n}-z^{-n}}{2i},b_{n-1}^{-1}\pi_{n-1}\rangle_{\mathbb{R}}.
\end{equation}

Throughout, as a convention, take $\sigma_{0}=1$, $\pi_{0}=0$ and
$\beta_{0}=0$ as well as $a_{0}=b_{0}=1$.

In deed, identifying the unit circle with the interval $[0,2\pi)$
via the map $\theta\rightarrow e^{i\theta}$, we get the
orthonormal trigonometric polynomials of the first class $\pi_{n}(\theta)$ and
$\sigma_{n}(\theta)$ for the over $\mathbb{R}$ linearly ordered trigonometric system
\begin{equation}
\{1, \sin\theta, \cos\theta, \ldots, \sin n\theta, \cos n\theta,
\ldots\}
\end{equation}
by the above process when $z=e^{i\theta},\,\theta\in [0,2\pi)$.

As noted in the introduction, by the uniqueness of solution of the the RHC (1.3), we have the following mutual representation theorem for OPUC and OTP.

\begin{thm}[\!\!\cite{dd08}]
Let $\mu$ be a nontrivial probability measure on the unit circle
$\partial \mathbb{D}$, $\{1, \pi_{n}, \sigma_{n}\}$ be
the unique system of the orthonormal Laurent polynomials of the first class on
the unit circle with respect to $\mu$, and $\{\Phi_{n}\}$ be the
unique system of the monic orthogonal polynomials on the unit circle
with respect to $\mu$. Then for any $z\in \mathbb{C}$ and $n\in
\mathbb{N}$,
\begin{equation}
\Phi_{2n-1}(z)=z^{n-1}[a_{n}\sigma_{n}(z)+(\beta_{n}+i)b_{n}\pi_{n}(z)]
\end{equation}
and
\begin{equation}
\kappa^{2}_{2n}\Phi^{*}_{2n}(z)=\frac{1}{2}z^{n}[a^{-1}_{n}(1+\beta_{n}i)\sigma_{n}(z)
-ib^{-1}_{n}\pi_{n}(z)],
\end{equation}
where $\kappa_{n}$ is the leading coefficient of the orthonormal
polynomial of order $n$ on the unit circle with respect to $\mu$, $\kappa_{n}=\|\Phi_{n}\|^{-1}_{\mathbb{C}}$, $\Phi_{n}^{*}$ is the reversed polynomial of $\Phi_{n}$, and $a_{n},
b_{n}, \beta_{n}$ are given in (2.11)-(2.13).
\end{thm}

Denote
$\Lambda_{n}=-\frac{1}{2}[a_{n}^{-2}(1+\beta_{n}^{2})+b_{n}^{-2}]i$,
by (2.17) and (2.18), we obtain
\begin{thm}
\begin{equation}
a_{n}\sigma_{n}(z)=-\frac{1}{2}z^{-n}[\Lambda_{n}^{-1}b_{n}^{-2}iz\Phi_{2n-1}(z)-(1-\beta_{n}i)\Phi_{2n}^{*}(z)]
\end{equation}
and
\begin{equation}
b_{n}\pi_{n}(z)=-\frac{1}{2}z^{-n}[\Lambda_{n}^{-1}a_{n}^{-2}(1+\beta_{n}i)z\Phi_{2n-1}(z)-i\Phi_{2n}^{*}(z)]
\end{equation}
for $n\in \mathbb{N}$ and $z\in \mathbb{C}\setminus\{0\}$.
\end{thm}

As some consequences, we have
\begin{thm} [\!\!\cite{dd08}]
\begin{equation}
\kappa_{2n}^{2}=\frac{1}{4}[a_{n}^{-2}(1+\beta_{n}^{2})+b_{n}^{-2}]
\end{equation}
for $n\in \mathbb{N}\cup\{0\}$.
\end{thm}

\begin{thm}
\begin{equation}
\alpha_{2n-1}=\frac{1}{4}\kappa_{2n}^{-2}[b_{n}^{-2}-a_{n}^{-2}(1-\beta_{n}^{2})]-\frac{1}{2}\kappa_{2n}^{-2}a_{n}^{-2}
\beta_{n}i
\end{equation}
and
\begin{equation}
\alpha_{2n-2}=\frac{1}{2}(\imath_{n}+\beta_{n-1}\varsigma_{n}-\zeta_{n})-\frac{i}{2}(\jmath_{n}-\imath_{n}\beta_{n-1}
+\varsigma_{n})
\end{equation}
for $n\in \mathbb{N}$.
\end{thm}
\begin{proof}
(2.22) is referred to \cite{dd08}. (2.23) follows from (2.11),
(2.12), (2.17) and the fact
$\alpha_{2n-2}=-\overline{\Phi_{2n-1}(0)}.$
\end{proof}

Since $\kappa_{n}^{2}/\kappa_{n+1}^{2}=1-|\alpha_{n}|^{2}$ for $n\in
\mathbb{N}\cup\{0\}$, by Theorem 2.3 and 2.4, we get

\begin{thm}
\begin{equation}
\kappa_{2n-1}^{2}=[a_{n}^{2}+b_{n}^{2}(1+\beta^{2}_{n})]^{-1}
\end{equation}
for $n\in \mathbb{N}$.
\end{thm}

Therefore, by (2.21) and (2.24), we obtain
\begin{thm}
\begin{equation}
\lim_{n\rightarrow
\infty}a_{n}b_{n}=\frac{1}{2}\exp\Big(\frac{1}{2\pi i}\int_{\partial
\mathbb{D}}\log w(\tau)\frac{d\tau}{\tau}\Big)
\end{equation}
and
\begin{equation}
\lim_{n\rightarrow
\infty}[a_{n}^{2}+b_{n}^{2}(1+\beta^{2}_{n})]=\exp\Big(\frac{1}{2\pi
i}\int_{\partial \mathbb{D}}\log w(\tau)\frac{d\tau}{\tau}\Big).
\end{equation}
\end{thm}
\begin{proof}
Since (see \cite{sze,sim1})
\begin{equation}
\lim_{n\rightarrow \infty}\kappa_{n}^{-2}=\exp\Big(\frac{1}{2\pi
i}\int_{\partial \mathbb{D}}\log w(\tau)\frac{d\tau}{\tau}\Big),
\end{equation}
then (2.26) follows from (2.24) whereas (2.25) follows from
\begin{equation}
\kappa_{2n-1}^{2}\kappa_{2n}^{2}=\frac{1}{4}a_{n}^{-2}b_{n}^{-2}
\end{equation}
which holds by (2.21) and (2.24).
\end{proof}

In addition, we also have
\begin{thm}
\begin{align}
&[a_{n}^{-2}(1+\beta_{n}^{2})+b_{n}^{-2}][a_{n+1}^{2}+b_{n+1}^{2}(1+\beta_{n+1}^{2})]+(\imath_{n+1}+\beta_{n}\varsigma_{n+1}-\zeta_{n+1})^{2}\nonumber\\
&+(\jmath_{n+1}-\imath_{n+1}\beta_{n}
+\varsigma_{n+1})^{2}=4
\end{align}
for $n\in \mathbb{N}\cup\{0\}$.
\end{thm}
\begin{proof}
It immediately follows from (2.21), (2.23) and (2.24) since
$\kappa_{2n}^{2}/\kappa_{2n+1}^{2}=1-|\alpha_{2n}|^{2}$ for $n\in
\mathbb{N}\cup\{0\}$.
\end{proof}

In the rest of this section, we give another identity on the coefficients $a_{n}, b_{n},\beta_{n}$ of OTP and $\alpha_{n}, \kappa_{n}$ of OPUC. The main idea will also be used in Section 5 below. To do so, we need the following simple facts.
\begin{lem}
Let $\Phi_{n}$ be the monic orthogonal polynomial on the unit circle of order $n$ with respect to $\mu$, and $\Phi^{*}_{n}$ be the reversed polynomial of $\Phi_{n}$, then
\begin{equation}
\langle 1,\Phi_{n}^{*}\rangle_{\mathbb{R}}=\int_{\partial \mathbb{D}}\Phi_{n}^{*}(\tau)d\mu(\tau)=\kappa_{n}^{-2}
\end{equation}
and
\begin{equation}
\langle 1,z\Phi_{n}\rangle_{\mathbb{R}}=\int_{\partial \mathbb{D}}\tau\Phi_{n}(\tau)d\mu(\tau)=\alpha_{n}^{-1}\Big(\kappa_{n}^{-2}-\kappa_{n+1}^{-2}\Big),
\end{equation}
where the Verblunsky coefficient $\alpha_{n}$ is restricted in $\mathbb{D}\setminus\{0\}$.
\end{lem}
\begin{proof}
By (2.6), we have
\begin{align}
\int_{\partial \mathbb{D}}\Phi_{n}^{*}(\tau)d\mu(\tau)=\int_{\partial \mathbb{D}}\tau^{n}\overline{\Phi_{n}(\tau)}d\mu(\tau)=\langle z^{n},\Phi_{n}\rangle_{\mathbb{C}}=\langle \Phi_{n},\Phi_{n}\rangle_{\mathbb{C}}.
\end{align}
Thus (2.30) holds on account of (2.4). For $\alpha_{n}\in\mathbb{D}\setminus\{0\}$, by Szeg\"o recurrence (2.7) (or see (4.12) below),
\begin{align}
\int_{\partial \mathbb{D}}\tau\Phi_{n}(\tau)d\mu(\tau)=\alpha_{n}^{-1}\Big[\int_{\partial \mathbb{D}}\Phi_{n}^{*}(\tau)d\mu(\tau)-\int_{\partial \mathbb{D}}\Phi_{n+1}^{*}(\tau)d\mu(\tau)\Big].
\end{align}
Therefore, (2.31) follows from (2.30).
\end{proof}

\begin{thm}
\begin{align}
&\alpha_{2n-1}\beta_{n}+\frac{1}{4}\Lambda_{n}^{-1}a_{n}^{-2}b_{n}^{-2}(1+\beta_{n}i)(1+\alpha_{2n-1})\kappa_{2n-1}^{-2}\nonumber
\\-&\frac{1}{4}\Big[\Lambda_{n}^{-1}a_{n}^{-2}b_{n}^{-2}(1+\beta_{n}i)+\alpha_{2n-1}b_{n}^{-2}i\Big]\kappa_{2n}^{-2}=0
\end{align}
for $n\in \mathbb{N}$.
\end{thm}
\begin{proof}
In the case of $\alpha_{2n-1}=0$, since $\kappa_{2n-1}=\kappa_{2n}$ by $\kappa_{2n-1}^{2}/\kappa_{2n}^{2}=1-|\alpha_{2n-1}|^{2}$, it is easy to get (2.34).
So in what follows, we always assume that $\alpha_{2n-1}\in \mathbb{D}\setminus\{0\}$.

By Theorem 2.2, Lemma 2.8 and Szeg\"o recurrence, we have
\begin{align}
\langle z^{n},b_{n}\pi_{n}\rangle_{\mathbb{R}}=&\langle z^{n},-\frac{1}{2}z^{-n}[\Lambda_{n}^{-1}a_{n}^{-2}(1+\beta_{n}i)z\Phi_{2n-1}-i\Phi_{2n}^{*}]\rangle_{\mathbb{R}}\nonumber\\
=&-\frac{1}{2}\Lambda_{n}^{-1}a_{n}^{-2}(1+\beta_{n}i)\langle 1,z\Phi_{2n-1}\rangle_{\mathbb{R}}+\frac{i}{2}\langle 1,\Phi_{2n}^{*}\rangle_{\mathbb{R}}\nonumber\\
=&-\frac{1}{2}\Lambda_{n}^{-1}a_{n}^{-2}(1+\beta_{n}i)\alpha_{2n-1}^{-1}\kappa_{2n-1}^{-2}+\Big[\frac{1}{2}\Lambda_{n}^{-1}a_{n}^{-2}(1+\beta_{n}i)\alpha_{2n-1}^{-1}+\frac{i}{2}\Big]\kappa_{2n}^{-2}\nonumber
\end{align}
and
\begin{align}
\langle z^{-n},b_{n}\pi_{n}\rangle_{\mathbb{R}}=&\langle z^{-n},-\frac{1}{2}z^{-n}[\Lambda_{n}^{-1}a_{n}^{-2}(1+\beta_{n}i)z\Phi_{2n-1}-i\Phi_{2n}^{*}]\rangle_{\mathbb{R}}\nonumber\\
=&-\frac{1}{2}\Lambda_{n}^{-1}a_{n}^{-2}(1+\beta_{n}i)\langle 1,z^{-(2n-1)}\Phi_{2n-1}\rangle_{\mathbb{R}}+\frac{i}{2}\langle 1,z^{-2n}\Phi_{2n}^{*}\rangle_{\mathbb{R}}\nonumber\\
=&-\frac{1}{2}\Lambda_{n}^{-1}a_{n}^{-2}(1+\beta_{n}i)\langle z^{(2n-1)},\Phi_{2n-1}\rangle_{\mathbb{C}}+\frac{i}{2}\overline{\langle 1,\Phi_{2n}\rangle}_{\mathbb{R}}\nonumber\\
=&-\frac{1}{2}\Lambda_{n}^{-1}a_{n}^{-2}(1+\beta_{n}i)\langle \Phi_{2n-1},\Phi_{2n-1}\rangle_{\mathbb{C}}\nonumber\\
=&-\frac{1}{2}\Lambda_{n}^{-1}a_{n}^{-2}(1+\beta_{n}i)\kappa_{2n-1}^{-2}.
\end{align}
Thus
\begin{align}
\beta_{n}=&\langle\frac{z^{n}+z^{-n}}{2},b_{n}^{-1}\pi_{n}\rangle_{\mathbb{R}}=b_{n}^{-2}\langle\frac{z^{n}+z^{-n}}{2},b_{n}\pi_{n}\rangle_{\mathbb{R}}\nonumber\\
=&-\frac{1}{4}\Lambda_{n}^{-1}a_{n}^{-2}b_{n}^{-2}(1+\beta_{n}i)\Big(\alpha_{2n-1}^{-1}+1\Big)\kappa_{2n-1}^{-2}\nonumber\\
&+\frac{1}{4}\Big[\Lambda_{n}^{-1}a_{n}^{-2}b_{n}^{-2}(1+\beta_{n}i)\alpha_{2n-1}^{-1}+b_{n}^{-2}i\Big]\kappa_{2n}^{-2}.
\end{align}
Multiplying by $\alpha_{2n-1}$ on two sides of (2.36), (2.34) immediately follows.
\end{proof}

\begin{rem}
Noting
\begin{equation}
\Lambda_{n}=-\frac{1}{2}[a_{n}^{-2}(1+\beta_{n}^{2})+b_{n}^{-2}]i=-2\kappa_{2n}^{2}i,
\end{equation}
we can also get (2.34) by directly invoking (2.21), (2.22) and (2.24) together.
\end{rem}

\section{Favard, Baxter, Geronimus,
Rakhmanov and Szeg\"o theorems}

In the present section, some theorems are obtained for orthogonal trigonometric polynomials, such
as Favard, Baxter, Geronimus, Rakhmanov theorems and so on, which play important roles in the theory
of OPUC
\cite{sim1,sim2}.

\subsection{Weak Favard Theorem} We begin with a weak Favard Theorem for OTP.
Favard theorem for OPRL is about
the orthogonality of a system of polynomials which satisfies a
three-term recurrence with appropriate coefficients \cite{sze,ma}.
Its OPUC version is well-known and also called Verblunsky theorem
\cite{sim1,enzg}, that is, if $\{\alpha_{n}^{(0)}\}_{n=0}^{\infty}$ is a
sequence of complex numbers in $\mathbb{D}$, then there exists a
unique measure $d\mu$ such that $\alpha_{n}(d\mu)=\alpha_{n}^{(0)}$,
where $\alpha_{n}(d\mu)$ are the associated Verblunsky coefficients
of $d\mu$.

For orthogonal trigonometric polynomials, we have the following Favard theorem in a weak form.
\begin{thm} Let
$\{(a_{n}^{(0)},b_{n}^{(0)},\beta_{n}^{(0)})\}_{n=0}^{\infty}$ with
$a_{0}^{(0)},b_{0}^{(0)}=1$ and $\beta_{0}^{(0)}=0$ be a system of
three-tuples of real numbers satisfying
\begin{align}
&[(a_{n}^{(0)})^{2}+(b_{n}^{(0)})^{2}(1+(\beta_{n}^{(0)})^{2})]
[(a_{n+1}^{(0)})^{2}+(b_{n+1}^{(0)})^{2}(1+(\beta_{n+1}^{(0)})^{2})]<4(a_{n}^{(0)})^{2}(b_{n}^{(0)})^{2}
\end{align}
with $a_{n}^{(0)},b_{n}^{(0)}>0$ for $n\in \mathbb{N}\cup\{0\}$,
then there exists a nontrivial probability measure $d\mu$ on
$\partial \mathbb{D}$ such that $a_{n}(d\mu)=a_{n}^{(0)}$,
$b_{n}(d\mu)=b_{n}^{(0)}$ and $\beta_{n}(d\mu)=\beta_{n}^{(0)}$,
where $a_{n}(d\mu),b_{n}(d\mu),\beta_{n}(d\mu)$ are associated
coefficients of $d\mu$ defined by (2.11)-(2.13).
\end{thm}

\begin{proof}
For $n\in \mathbb{N}\cup\{0\}$, define \begin{equation}
\kappa_{2n}^{(0)}=\frac{1}{2}\Big[(a_{n}^{(0)})^{-2}\big(1+(\beta_{n}^{(0)})^{2}\big)+(b_{n}^{(0)})^{-2}\Big]^{\frac{1}{2}}
\end{equation}
and
\begin{equation}
\kappa_{2n+1}^{(0)}=\Big[(a_{n+1}^{(0)})^{2}+(b_{n+1}^{(0)})^{2}\big(1+(\beta_{n+1}^{(0)})^{2}\big)\Big]^{-\frac{1}{2}}.
\end{equation}

Let
\begin{equation}
\alpha_{2n-1}^{(0)}=\frac{1}{4}(\kappa_{2n}^{(0)})^{-2}\Big[(b_{n}^{(0)})^{-2}-(a_{n}^{(0)})^{-2}\big(1-(\beta_{n}^{(0)})^{2}\big)\Big]
-\frac{1}{2}(\kappa_{2n}^{(0)})^{-2}(a_{n}^{(0)})^{-2}
(\beta_{n}^{(0)})i,
\end{equation}
then $\alpha_{2n-1}^{(0)}\in \mathbb{D}$ since
\begin{equation}
\Big|\alpha_{2n-1}^{(0)}\Big|^{2}=\frac{(\kappa_{2n}^{(0)})^{4}-\frac{1}{4}(a_{n}^{(0)})^{-2}(b_{n}^{(0)})^{-2}}
{(\kappa_{2n}^{(0)})^{4}}=1-\frac{(\kappa_{2n-1}^{(0)})^{2}}{(\kappa_{2n}^{(0)})^{2}}
\end{equation}
and $a_{n}^{(0)},b_{n}^{(0)}>0$. Note that (3.1) is equivalent to
\begin{equation}
\frac{\kappa_{2n}^{(0)}}{\kappa_{2n+1}^{(0)}}<1.
\end{equation}

Arbitrarily choose a sequence $\{\alpha_{2n}^{(0)}\}_{n=0}^{\infty}$
such that
\begin{equation}
\Big|\alpha_{2n}^{(0)}\Big|=\sqrt{1-(\kappa_{2n}^{(0)})^{2}\big/(\kappa_{2n+1}^{(0)})^{2}}
\end{equation}
and fix it, then $\alpha_{2n}^{(0)}\in\mathbb{D}$ for $n\in
\mathbb{N}\cup\{0\}$ by (3.6).

Therefore, for this fixed sequence
$\{\alpha_{n}^{(0)}\}_{n=0}^{\infty}$, by Verblunsky theorem, there
exists a unique nontrivial probability measure $d\mu$ on $\partial
\mathbb{D}$ such that
\begin{equation}
\alpha_{n}(d\mu)=\alpha_{n}^{(0)}
\end{equation} for $n\in
\mathbb{N}\cup\{0\}$. Then for $n\in \mathbb{N}\cup\{0\}$,
\begin{equation}
\kappa_{n}(d\mu)=\kappa_{n}^{(0)}
\end{equation}
since
$\kappa_{n}(d\mu)=\prod_{j=0}^{n-1}(1-|\alpha_{j}(d\mu)|^{2})^{-\frac{1}{2}}$
(see \cite{sim1}).

Suppose that $\{\Phi_{n}(d\mu,z)\}_{n=0}^{\infty}$ is the sequence
of monic orthogonal polynomials on the unit circle with respect to
$d\mu$, set
\begin{equation}
\Sigma_{n}(z)=-\frac{1}{2}z^{-n}[(\Lambda_{n}^{(0)})^{-1}(b_{n}^{(0)})^{-2}iz\Phi_{2n-1}(d\mu,z)
-(1-\beta_{n}^{(0)}i)\Phi_{2n}^{*}(d\mu,z)]
\end{equation}
and
\begin{equation}
\Pi_{n}(z)=-\frac{1}{2}z^{-n}[(\Lambda_{n}^{(0)})^{-1}(a_{n}^{(0)})^{-2}(1+\beta_{n}^{(0)}i)z\Phi_{2n-1}(d\mu,z)
-i\Phi_{2n}^{*}(d\mu,z)]
\end{equation}
for $n\in \mathbb{N}$ and $z\in \mathbb{C}\setminus\{0\}$, where
$\Lambda_{n}^{(0)}=-\frac{1}{2}\Big[(a_{n}^{(0)})^{-2}\big(1+(\beta_{n}^{(0)})^{2}\big)+(b_{n}^{(0)})^{-2}\Big]i$.
Obviously,
\begin{equation}
\Lambda_{n}^{(0)}=-2(\kappa_{2n}^{(0)})^{2}i.
\end{equation}
By Szeg\"o recurrence and (3.8),
\begin{equation}
z\Phi_{2n-1}(d\mu,z)=\Phi_{2n}(d\mu,z)+\overline{\alpha^{(0)}_{2n-1}}\Phi^{*}_{2n-1}(d\mu,z).
\end{equation}
Hence by the orthogonality of $\Phi_{n}(d\mu, z)$ and
$\Phi_{n}^{*}(d\mu, z)$, we get
\begin{equation}
\langle z^{\pm j}, \Sigma_{n}\rangle_{\mathbb{R}}=\langle z^{\pm j},
\Pi_{n}\rangle_{\mathbb{R}}=0,\,\,\,\,j=0,1,\ldots,n-1.
\end{equation}
Moreover,
\begin{equation}
\langle z^{n},
\Sigma_{n}\rangle_{\mathbb{R}}=(a_{n}^{(0)})^{2}\overline{\alpha^{(0)}_{2n-1}}+\frac{1}{2}(\kappa_{2n}^{(0)})^{-2}(1-\beta_{n}^{(0)}i),
\end{equation}
\begin{equation}
\langle z^{-n}, \Sigma_{n}\rangle_{\mathbb{R}}=(a_{n}^{(0)})^{2},
\end{equation}
\begin{equation}
\langle z^{n},
\Pi_{n}\rangle_{\mathbb{R}}=(b_{n}^{(0)})^{2}(\beta_{n}^{(0)}-i)\overline{\alpha^{(0)}_{2n-1}}+\frac{1}{2}(\kappa_{2n}^{(0)})^{-2}i,
\end{equation}
and
\begin{equation}
\langle z^{-n},
\Pi_{n}\rangle_{\mathbb{R}}=(b_{n}^{(0)})^{2}(\beta_{n}^{(0)}-i)
\end{equation}
follow from (3.9), (3.12) and the fact
$||\Phi_{n}(d\mu)||_{\mathbb{R}}^{2}=||\Phi_{n}^{*}(d\mu)||_{\mathbb{R}}^{2}=[\kappa_{n}(d\mu)]^{^{-2}}$
as well as
$(\kappa_{2n-1}^{(0)})^{2}(\kappa_{2n}^{(0)})^{2}=\frac{1}{4}(a_{n}^{(0)})^{-2}(b_{n}^{(0)})^{-2}$.
By (3.4),
\begin{equation}
\overline{\alpha^{(0)}_{2n-1}}-1=-\frac{1}{2}(\kappa_{2n}^{(0)})^{-2}(a_{n}^{(0)})^{-2}(1-\beta_{n}^{(0)}i)
\end{equation}
and
\begin{equation}
\overline{\alpha^{(0)}_{2n-1}}+1=\frac{1}{2}(\kappa_{2n}^{(0)})^{-2}\Big[(a_{n}^{(0)})^{-2}(\beta_{n}^{(0)})^{2}
+(b_{n}^{(0)})^{-2}\Big]+\frac{1}{2}(\kappa_{2n}^{(0)})^{-2}(a_{n}^{(0)})^{-2}\beta_{n}^{(0)}i.
\end{equation}
So
\begin{equation}
\langle \frac{z^{n}+z^{-n}}{2},
\Sigma_{n}\rangle_{\mathbb{R}}=(a_{n}^{(0)})^{2},\,\,\,\langle
\frac{z^{n}-z^{-n}}{2i},
\Pi_{n}\rangle_{\mathbb{R}}=(b_{n}^{(0)})^{2}
\end{equation}
and
\begin{equation}
\langle \frac{z^{n}-z^{-n}}{2i}, \Sigma_{n}\rangle_{\mathbb{R}}=0
\end{equation}
as well as
\begin{equation}
\langle \frac{z^{n}+z^{-n}}{2},
\Pi_{n}\rangle_{\mathbb{R}}=(b_{n}^{(0)})^{2}\beta_{n}^{(0)}.
\end{equation}

In addition, it is easy to check that the
coefficients of $z^{n}$ and $z^{-n}$ in $\Pi_{n}(z)$ are
respectively $\frac{1}{2i}$ and $-\frac{1}{2i}$ whereas both of ones
in $\Sigma_{n}(z)-\beta_{n}^{(0)}\Pi_{n}(z)$ are $\frac{1}{2}$. Noting
(3.14) and (3.22), this fact means that $\Sigma_{n}(z)$ and
$\Pi_{n}(z)$ are just the ``monic" orthogonal Laurent
polynomials of the first class on the unit circle with respect to $d\mu$, i.e.
\begin{equation}
\Sigma_{n}(z)=a_{n}(d\mu)\sigma_{n}(d\mu,z)\,\,\,\,\,\text{and}\,\,\,\,\,\Pi_{n}(z)=b_{n}(d\mu)\pi_{n}(d\mu,z).
\end{equation}
Since
\begin{equation}
\langle a_{n}(d\mu)\sigma_{n}(d\mu),
a_{n}(d\mu)\sigma_{n}(d\mu)\rangle_{\mathbb{R}}=a_{n}^{2}(d\mu),
\end{equation}
\begin{equation}
\langle b_{n}(d\mu)\pi_{n}(d\mu),
b_{n}(d\mu)\pi_{n}(d\mu)\rangle_{\mathbb{R}}=b_{n}^{2}(d\mu)
\end{equation}
and
\begin{equation}
\langle \frac{z^{n}+z^{-n}}{2},
b_{n}(d\mu)\pi_{n}(d\mu)\rangle_{\mathbb{R}}=b_{n}^{2}(d\mu)\beta_{n}(d\mu),
\end{equation}
therefore, by (3.21) and (3.23),
\begin{equation}
a_{n}(d\mu)=a_{n}^{(0)},\,\,\,b_{n}(d\mu)=b_{n}^{(0)},\,\,\,\beta_{n}(d\mu)=\beta_{n}^{(0)}.
\end{equation}
\end{proof}

\begin{rem}
Only for the sequence of three-tuples
$(a_{n}^{(0)},b_{n}^{(0)},\beta_{n}^{(0)})$ fulfilling (3.1), to get
(3.28), the measure $d\mu$ is not unique since the sequence can
definitely determine Verblunsky coefficients with odd subscript but
ones with even subscript from the above proof.

For $n\in \mathbb{N}$, set
\begin{equation}
\imath_{n}(d\mu)=\langle\frac{z^{n}+z^{-n}}{2},(a_{n-1}^{(0)})^{-1}\sigma_{n-1}(d\mu)\rangle_{\mathbb{R}},
\end{equation}
\begin{equation}
\jmath_{n}(d\mu)=\langle\frac{z^{n}+z^{-n}}{2},(b_{n-1}^{(0)})^{-1}\pi_{n-1}(d\mu)\rangle_{\mathbb{R}},
\end{equation}
\begin{equation}
\varsigma_{n}(d\mu)=\langle\frac{z^{n}-z^{-n}}{2i},(a_{n-1}^{(0)})^{-1}\sigma_{n-1}(d\mu)\rangle_{\mathbb{R}},
\end{equation}
and
\begin{equation}
\zeta_{n}(d\mu)=\langle\frac{z^{n}-z^{-n}}{2i},(b_{n-1}^{(0)})^{-1}\pi_{n-1}(d\mu)\rangle_{\mathbb{R}}.
\end{equation}
Then the measure $d\mu$ is unique for the sequence of seven-tuples
\begin{equation}
(a_{n}^{(0)},b_{n}^{(0)},\beta_{n}^{(0)},\imath_{n}(d\mu),\jmath_{n}(d\mu),\varsigma_{n}(d\mu),\zeta_{n}(d\mu))
\end{equation}
satisfying (3.1) by Theorem 2.4 and Verblunsky theorem. Since $d\mu$
is partly dependent on $(a_{n}^{(0)},b_{n}^{(0)},\beta_{n}^{(0)})$ and
$\imath_{n}(d\mu),\jmath_{n}(d\mu),\varsigma_{n}(d\mu),\zeta_{n}(d\mu)$
are dependent on $d\mu$, $a_{n}^{(0)}$ and $b_{n}^{(0)}$, then the
sequence of seven-tuples (3.33) satisfying (3.1) is partly dependent on the
sequence of three-tuples $(a_{n}^{(0)},b_{n}^{(0)},\beta_{n}^{(0)})$
fulfilling (3.1). Considering the uniqueness of $d\mu$ for the
sequence of (3.33) with (3.1), we call that $d\mu$ is selectively
unique for the sequence
$\{(a_{n}^{(0)},b_{n}^{(0)},\beta_{n}^{(0)})\}_{n=0}^{\infty}$
satisfying (3.1) and $a_{n}^{(0)},b_{n}^{(0)}>0$ as well as
$a_{0}^{(0)},b_{0}^{(0)}=1$ and $\beta_{0}^{(0)}=0$. In Section 5, we will give a strong Favard theorem which in detail illuminates the relation on the uniqueness of $d\mu$ and a sequence of seven-tuples, $\{(a_{n}^{(0)}, b_{n}^{(0)},\beta_{n}^{(0)},\imath_{n}^{(0)},\jmath_{n}^{(0)},\varsigma_{n}^{(0)},\zeta_{n}^{(0)})\}$, with some additional properties.
\end{rem}

Similarly, by Theorems 2.3, 2.4, 2.7 and using the corresponding theorems for OPUC, we also have Baxter, Geronimus,
Rakhmanov, Szeg\"o and the strong Szeg\"o theorems for OTP in what follows.
\subsection{Baxter Theorem}
Let
\begin{equation}
c_{n}=\int_{\partial \mathbb{D}}\overline{\tau}^{n}d\mu(\tau),
\,\,\,n\in \mathbb{N}\cup\{0\}
\end{equation}
be moments of $\mu$, Baxter theorem for OPUC states that
$\sum_{n=0}^{\infty}|\alpha_{n}|<0$ if and only if
$\sum_{n=0}^{\infty}|c_{n}|<0$ and
$d\mu(\tau)=w(\tau)\frac{d\tau}{2\pi i\tau}$ with $w(\tau)$
continuous and $\min_{\tau\in\partial \mathbb{D}}w(\tau)>0$.

For orthogonal trigonometric polynomials, we have Baxter theorem as follows.
\begin{thm}
Let $\mu$ be a nontrivial probability measure on $\partial
\mathbb{D}$, $a_{n}, b_{n}, \beta_{n}$ be the associated
coefficients given in (2.11)-(2.13) and $c_{n}$ be the moments of
$\mu$ defined by (3.34), then
\begin{align}
&\sum_{n=0}^{\infty}\sqrt{1-\frac{1}{4}[a_{n}^{-2}(1+\beta_{n}^{2})+b_{n}^{-2}][a_{n+1}^{2}+b_{n+1}^{2}(1+\beta_{n+1}^{2})]}
\nonumber\\
+&\sum_{n=0}^{\infty}\sqrt{\frac{a_{n}^{4}+b_{n}^{4}(1+\beta_{n}^{2})^{2}+2a_{n}^{2}b_{n}^{2}(\beta_{n}^{2}-1)}
{a_{n}^{4}+b_{n}^{4}(1+\beta_{n}^{2})^{2}+2a_{n}^{2}b_{n}^{2}(\beta_{n}^{2}+1)}}<\infty
\end{align}
is equivalent to $\sum_{n=0}^{\infty}|c_{n}|<0$ and
$d\mu(\tau)=w(\tau)\frac{d\tau}{2\pi i\tau}$ with $w(\tau)$
continuous and $\min_{\tau\in\partial \mathbb{D}}w(\tau)>0$.
\end{thm}

\subsection{Geronimus Theorem} To discuss Geronimus theorem, it is
necessary to introduce some basic notions of Schur algorithm (see
\cite{sim1}).

An analytic function $F$ on $\mathbb{D}$ is called a Carath\'eodory
function if and only if $F(0)=1$ and $\Re F(z)>0$ on $\mathbb{D}$.
An analytic function $f$ on $\mathbb{D}$ is called a Schur function
if and only if $\sup_{z\in \mathbb{D}}|f(z)|<1$. Let
\begin{equation}
F(z)=\int_{\partial\mathbb{D}}\frac{\tau+z}{\tau-z}d\mu(\tau)
\end{equation}
be an associated Carath\'eodory function of $\mu$, then
\begin{equation}
f(z)=\frac{1}{z}\frac{F(z)-1}{F(z)+1}
\end{equation}
is a Schur function related to $\mu$.

Starting with a Schur function $f_{0}$, Schur algorithm actually
provides an approach to continuously map one Schur function to
another by a series of transforms of the form
\begin{equation}
\begin{cases}
f_{n+1}(z)=\displaystyle\frac{1}{z}\frac{f_{n}(z)-\gamma_{n}}{1-\overline{\gamma}_{n}f_{n}(z)},\\[4mm]
\gamma_{n}=f_{n}(0).
\end{cases}
\end{equation}
$f_{n}$ are called Schur iterates and $\gamma_{n}$ are called Schur
parameters associated to $f_{0}$. Due to Schur, it is well known
that there is a one to one correspondence between the set of Schur
functions which are not finite Blaschke products and the set of
sequences of $\{\gamma_{n}\}_{n=0}^{\infty}$ in $\mathbb{D}$.
Geronimus theorem for OPUC asserts that if $\mu$ is a nontrivial
probability measure on $\partial \mathbb{D}$, the Schur parameters
$\{\gamma_{n}\}_{n=0}^{\infty}$ associated to $f_{0}$ related to
$\mu$ defined by (3.36) and (3.37) are identical to the Verblunsky
coefficients $\{\alpha_{n}\}_{n=0}^{\infty}$.

For orthogonal trigonometric polynomials, we have Geronimus theorem as follows.
\begin{thm}
Let $\mu$ be a nontrivial probability measure on $\partial
\mathbb{D}$, if $\gamma_{n}$ are Schur parameters and $a_{n}$,
$b_{n}$, $\beta_{n}$, $\imath_{n}$, $\jmath_{n}$, $\varsigma_{n}$,
$\zeta_{n}$ are coefficients associated to $\mu$ defined by
(2.11)-(2.15), then
\begin{equation}
\gamma_{2n-1}=\frac{a_{n}^{2}-b_{n}^{2}(1-\beta_{n}^{2})}{a_{n}^{2}+b_{n}^{2}(1+\beta_{n}^{2})}
-\frac{2b_{n}^{2}\beta_{n}}{a_{n}^{2}+b_{n}^{2}(1+\beta_{n}^{2})}i
\end{equation}
and
\begin{equation}
\gamma_{2n-2}=\frac{1}{2}(\imath_{n}+\beta_{n-1}\varsigma_{n}-\zeta_{n})-\frac{i}{2}(\jmath_{n}-\imath_{n}\beta_{n-1}
+\varsigma_{n})
\end{equation}
for $n\in \mathbb{N}$.
\end{thm}

\subsection{Rakhmanov Theorem and Szeg\"o Theorem}
Let $d\mu$ have the decomposition form (2.1),
$\{\alpha_{n}\}_{n=0}^{\infty}$ be the Verblunsky coefficients of
$\mu$, Rakhmanov theorem for OPUC states that if $w(\tau)>0$ for
a.e. $\tau\in
\partial\mathbb{D}$, then $\lim_{n\rightarrow\infty}|\alpha_{n}|=0$.
Its OTP version is as follows.
\begin{thm}
Let $\mu$ be a nontrivial probability measure on $\partial
\mathbb{D}$ with the decomposition form (2.1), $a_{n}, b_{n},
\beta_{n}$ be the associated coefficients of $\mu$ given in
(2.11)-(2.13). If $w(\tau)>0$ for a.e. $\tau\in \partial\mathbb{D}$,
then
\begin{equation}
\lim_{n\rightarrow\infty}\frac{a_{n}^{4}+b_{n}^{4}(1+\beta_{n}^{2})^{2}+2a_{n}^{2}b_{n}^{2}(\beta_{n}^{2}-1)}
{a_{n}^{4}+b_{n}^{4}(1+\beta_{n}^{2})^{2}+2a_{n}^{2}b_{n}^{2}(\beta_{n}^{2}+1)}=0
\end{equation}
and
\begin{equation}
\lim_{n\rightarrow\infty}\frac{1}{4}[a_{n}^{-2}(1+\beta_{n}^{2})+b_{n}^{-2}][a_{n+1}^{2}+b_{n+1}^{2}(1+\beta_{n+1}^{2})]=1.
\end{equation}
\end{thm}

Szeg\"o theorem for OPUC shows that
\begin{equation}
\prod_{n=0}^{\infty}(1-|\alpha_{n}|^{2})=\exp\Big(\frac{1}{2\pi
i}\int_{\partial \mathbb{D}}\log w(\tau)\frac{d\tau}{\tau}\Big).
\end{equation}

 Especially,
 \begin{equation}
 \sum_{n=0}^{\infty}|\alpha_{n}|^{2}<\infty\Longleftrightarrow
\frac{1}{2\pi i}\int_{\partial \mathbb{D}}\log
w(\tau)\frac{d\tau}{\tau}>-\infty.
 \end{equation}

Its analog for OTP is
\begin{thm}
Let $\mu$ be a nontrivial probability measure on $\partial
\mathbb{D}$ with the decomposition form (2.1), $a_{n}, b_{n},
\beta_{n}$ be the associated coefficients of $\mu$ given in
(2.11)-(2.13). Then
\begin{equation}
\prod_{n=0}^{\infty}\frac{a_{n+1}^{2}+b_{n+1}^{2}(1+\beta_{n+1}^{2})}{a_{n}^{2}+b_{n}^{2}(1+\beta_{n}^{2})}=
\exp\Big(\frac{1}{2\pi i}\int_{\partial \mathbb{D}}\log
w(\tau)\frac{d\tau}{\tau}\Big).
\end{equation}
In particular,
\begin{align}
&\sum_{n=0}^{\infty}\left\{1-\frac{1}{4}[a_{n}^{-2}(1+\beta_{n}^{2})+b_{n}^{-2}][a_{n+1}^{2}+b_{n+1}^{2}(1+\beta_{n+1}^{2})]\right\}
\nonumber\\
+&\sum_{n=0}^{\infty}\frac{a_{n}^{4}+b_{n}^{4}(1+\beta_{n}^{2})^{2}+2a_{n}^{2}b_{n}^{2}(\beta_{n}^{2}-1)}
{a_{n}^{4}+b_{n}^{4}(1+\beta_{n}^{2})^{2}+2a_{n}^{2}b_{n}^{2}(\beta_{n}^{2}+1)}<\infty
\end{align}
is equivalent to $\displaystyle\frac{1}{2\pi i}\int_{\partial
\mathbb{D}}\log w(\tau)\frac{d\tau}{\tau}>-\infty$.
\end{thm}

\subsection{The Strong Szeg\"o Theorem}Let $d\mu$ have the decomposition form
(2.1) satisfying the Szeg\"o condition
\begin{equation}
\frac{1}{2\pi i}\int_{\partial \mathbb{D}}\log
w(\tau)\frac{d\tau}{\tau}>-\infty,
\end{equation}
it is accustomed to introduce the Szeg\"o function as follows
\begin{equation}
D(z)=\exp\Big(\frac{1}{4\pi i}\int_{\partial
\mathbb{D}}\frac{\tau+z}{\tau-z}\log w(\tau)\frac{d\tau}{\tau}\Big).
\end{equation}
It is easy to get that $D(z)$ is analytic and nonvanishing in
$\mathbb{D}$, lies in the Hardy space $H^{2}(\mathbb{D})$ and
$\lim_{r\uparrow1}D(r\tau)=D(\tau)$ for a.e. $\tau\in\partial
\mathbb{D}$ as well as $|D(\tau)|^{2}=w(\tau)$. Let
\begin{equation}
D(z)=\exp\Big(\frac{1}{2}\hat{L}_{0}+\sum_{n=1}^{\infty}\hat{L}_{n}z^{n}\Big),\,\,\,z\in
\mathbb{D}.
\end{equation}
Due to Ibragimov, the sharpest form of the strong Szeg\"o theorem
for OPUC (see \cite{sim1}) says that
\begin{equation}
\sum_{n=0}^{\infty}n|\alpha_{n}|^{2}<\infty\Longleftrightarrow
d\mu_{s}=0 \,\,\,\text{and}
\,\,\,\sum_{n=0}^{\infty}n|\hat{L}_{n}|^{2}<\infty.
\end{equation}

The corresponding result for OTP is stated in the following theorem.
\begin{thm}
Let $\mu$ be a nontrivial probability measure on $\partial
\mathbb{D}$ with the decomposition form (2.1) satisfying the Szeg\"o
cindition (3.47), $a_{n}, b_{n}, \beta_{n}$ be the associated
coefficients of $\mu$ given in (2.11)-(2.13), and
$\{\hat{L}_{n}\}_{n=0}^{\infty}$ be the Taylor coefficients of the logarithm of
Szeg\"o function $D(z)$ at $z=0$ which are defined by (3.48) and
(3.49). Then
\begin{align}
&\sum_{n=0}^{\infty}2n\left\{1-\frac{1}{4}[a_{n}^{-2}(1+\beta_{n}^{2})+b_{n}^{-2}][a_{n+1}^{2}+b_{n+1}^{2}(1+\beta_{n+1}^{2})]\right\}
\nonumber\\
+&\sum_{n=0}^{\infty}(2n-1)\left\{\frac{a_{n}^{4}+b_{n}^{4}(1+\beta_{n}^{2})^{2}+2a_{n}^{2}b_{n}^{2}(\beta_{n}^{2}-1)}
{a_{n}^{4}+b_{n}^{4}(1+\beta_{n}^{2})^{2}+2a_{n}^{2}b_{n}^{2}(\beta_{n}^{2}+1)}\right\}<\infty
\end{align}
is equivalent to $d\mu_{s}=0$ and
$\sum_{n=0}^{\infty}n|\hat{L}_{n}|^{2}<\infty$.
\end{thm}

In the above, by the mutual representation theorem for OTP and OPUC,
we obtain some classical theorems for orthogonal trigonometric
polynomials corresponding to the ones for orthogonal polynomials on the
unit circle. In fact, by this theorem, we can
obtain much more results for orthogonal trigonometric
polynomials. For example, the important and useful Bernstein-Szeg\"o
measure can be expressed in terms of orthogonal trigonometric
polynomials as follows
\begin{equation*}
d\mu_{n}=
\begin{cases}
\displaystyle\frac{a_{m}^{2}+b_{m}^{2}(1+\beta^{2}_{m})}{|a_{m}\sigma_{m}(\theta)+(\beta_{m}+i)b_{m}\pi_{m}(\theta)|^{2}}
\frac{d\theta}{2\pi},\,\,\,n=2m-1,\\[3mm]
\displaystyle\frac{a_{m}^{2}b_{m}^{2}}{a_{m}^{2}+b_{m}^{2}(1+\beta^{2}_{m})}\frac{1}{|a_{m}^{-1}(\beta_{m}-i)\sigma_{m}(\theta)
-b_{m}^{-1}\pi_{m}(\theta)|^{2}} \frac{d\theta}{2\pi},\,\,\,n=2m.
\end{cases}
\end{equation*}

\section{Identities from Riemann-Hilbert Characterizations}

In this section, by applying the corresponding RH characterizations, we obtain some identities for OPUC and OTP including Szeg\"o recursions for OPUC, four-term recurrences for OTP and some new identities on Cauchy integrals for OPUC and OTP as well as Hilbert transforms for OTP. Let $H(\partial \mathbb{D})$ denote the set of all complex-valued and H\"older continuous functions defined on $\partial \mathbb{D}$. For simplicity, we always assume that the weight function $w\in H(\partial \mathbb{D})$ in what follows.

\subsection{The case of OPUC} The RH characterization  for OPUC is uniquely solvable as follows
\begin{thm}[\!\!\cite{bdj,dd06}]
The RHP (1.2) has a unique solution given by
\begin{equation}
Y(z)=\left(
  \begin{array}{cc}
    \Phi_{n}(z) & C[\tau^{-n}\Phi_{n}w](z) \\
    -\kappa_{n-1}^{2}\Phi_{n-1}^{*}(z) & -\kappa_{n-1}^{2}C[\tau^{-n}\Phi_{n-1}^{*}w](z) \\
  \end{array}
\right),
\end{equation}
where $\Phi_{n}$ is the monic orthogonal polynomials on the unit circle of order $n$ with respect to the weight $w$, $\Phi_{n-1}^{*}$ is the reversed polynomial of $\Phi_{n-1}$, $\kappa_{n-1}$ is given as in (2.4), and $C$ is the Cauchy integral operator.
\end{thm}

Consider the Schwarz reflection of $Y$ defined by
\begin{equation}
Y_{1}(z)=\overline{Y\left(\frac{1}{\overline{z}}\right)},\,\,\,\,z\in \mathbb{C}\setminus\partial \mathbb{D},
\end{equation}
then $Y_{1}^{+}(t)=\overline{Y^{-}(t)}$  and $Y_{1}^{-}(t)=\overline{Y^{+}(t)}$ for $t\in \partial \mathbb{D}$. Therefore, by using the boundary condition in RHP (1.2), we have
\begin{equation}
Y_{1}^{+}(t)=Y_{1}^{-}(t)\left(
                   \begin{array}{cc}
                     1 & -t^{n}w(t) \\
                     0 & 1 \\
                   \end{array}
                 \right),\,\,\, t\in \partial \mathbb{D}.
\end{equation}
By a direct evaluation,
\begin{equation}
\lim_{z\rightarrow \infty}Y_{1}(z)=\left(
                   \begin{array}{cc}
                     -\alpha_{n-1} & \kappa_{n}^{-2} \\
                     -\kappa_{n-1}^{2} & -\overline{\alpha}_{n-1} \\
                   \end{array}
                 \right).
\end{equation}
Moreover, by the growth condition at $\infty$ in RHP (1.2), we have
\begin{equation}
\lim_{z\rightarrow 0}Y_{1}(z)\left(
               \begin{array}{cc}
                     z^{n} & 0 \\
                     0 & z^{-n} \\
                   \end{array}
                 \right)=I.
\end{equation}
Let
\begin{align}
Y_{2}(z)=
\left(
               \begin{array}{cc}
                     -\overline{\alpha}_{n-1} & -\kappa_{n}^{-2} \\
                     -\kappa_{n-1}^{2} & \alpha_{n-1} \\
                   \end{array}
                 \right)Y_{1}(z)\left(
               \begin{array}{cc}
                     z^{n} & 0 \\
                     0 & -z^{-n} \\
                   \end{array}
                 \right).
\end{align}
Noting
\begin{equation}
1-|\alpha_{n-1}|^{2}=\left(\frac{\kappa_{n-1}}{\kappa_{n}}\right)^{2},
\end{equation}
by (4.4), (4.6) and simple calculations,
\begin{equation}
\lim_{z\rightarrow \infty}Y_{2}(z)\left(
               \begin{array}{cc}
                     z^{-n} & 0 \\
                     0 & z^{n} \\
                   \end{array}
                 \right)=I.
\end{equation}
By (4.3) and (4.8), $Y_{2}$ satisfies the RH characterization (1.2) for OPUC, viz.
\begin{equation} (\mbox{RHP for $Y_{2}$})\,\,
\begin{cases}
Y_{2}\,\, \mbox{is analytic in}\,\,\mathbb{C}\setminus \partial \mathbb{D},\vspace{2mm}\\
Y_{2}^{+}(t)=Y_{2}^{-}(t)\left(
                   \begin{array}{cc}
                     1 & t^{-n}w(t) \\
                     0 & 1 \\
                   \end{array}
                 \right)
                 \,\,\mbox{for} \,\,t\in \partial \mathbb{D},\vspace{2mm}\\
Y_{2}(z)=\left(I+O(\frac{1}{z})\right)\left(
                                    \begin{array}{cc}
                                      z^{n} & 0 \\
                                      0 & z^{-n} \\
                                    \end{array}
                                  \right)
                                  \,\,\mbox{as}\,\, z\rightarrow \infty.
\end{cases}
\end{equation}

By the uniqueness, $Y_{2}(z)=Y(z)$ for $z\in \mathbb{C}\setminus \partial \mathbb{D}$. Namely,

\begin{equation}
Y(z)=\left(
               \begin{array}{cc}
                     -\overline{\alpha}_{n-1} & -\kappa_{n}^{-2} \\
                     -\kappa_{n-1}^{2} & \alpha_{n-1} \\
                   \end{array}
                 \right)\overline{Y\left(\frac{1}{\overline{z}}\right)}\left(
               \begin{array}{cc}
                     z^{n} & 0 \\
                     0 &  -z^{-n}\\
                   \end{array}
                 \right),\,\,z\in \mathbb{C}\setminus\partial\mathbb{D}.
\end{equation}

By the above arguments, we have
\begin{thm}
Let $\Phi_{n}$, $\Phi_{n-1}^{*}$, $\alpha_{n-1}$, $\kappa_{n}$ be as above, then
\begin{enumerate}
  \item [(A)] The identities
  \begin{equation}
  \Phi_{n}(z)=-\overline{\alpha}_{n-1}\Phi_{n}^{*}(z)+\frac{\kappa_{n-1}^{2}}{\kappa_{n}^{2}}z\Phi_{n-1}(z)
  \end{equation}
  and
  \begin{equation}
  \Phi_{n}^{*}(z)=\Phi_{n-1}^{*}(z)-\alpha_{n-1}z\Phi_{n-1}(z)
  \end{equation}
  hold for $z\in \mathbb{C}$;
  \item [(B)] The identities
  \begin{align}
  C[\Phi_{n}w](z)=&z^{n}\Big[\overline{\alpha}_{n-1}\overline{C[\Phi_{n}w]\left(\frac{1}{\overline{z}}\right)}
  -\frac{\kappa_{n-1}^{2}}{\kappa_{n}^{2}}\overline{C[\Phi_{n-1}^{*}w]\left(\frac{1}{\overline{z}}\right)}+\frac{1}{\kappa_{n}^{2}}\Big]
  \end{align}
  and
  \begin{align}
  C[\Phi_{n-1}^{*}w](z)=&-z^{n}\Big[\overline{C[\Phi_{n}w]\left(\frac{1}{\overline{z}}\right)}
 +\alpha_{n-1}\overline{C[\Phi_{n-1}^{*}w]\left(\frac{1}{\overline{z}}\right)}-\frac{1+\alpha_{n-1}}{\kappa_{n-1}^{2}}\Big]
  \end{align}
  hold for $z\in \mathbb{C}\setminus\partial\mathbb{D}$.
\end{enumerate}
\end{thm}

\begin{proof}
(4.11) and (4.12) are obtained by identifying the 11 and 21 entries in left hand side with the ones in right hand side of (4.10).

By the orthogonality, it is easy to get that (see \cite{dd06})
\begin{equation}
C[\tau^{-n}\Phi_{n}w](z)=z^{-n}C[\Phi_{n}w](z)
\end{equation}
and
\begin{equation}
-\kappa_{n-1}^{2}C[\tau^{-n}\Phi_{n-1}^{*}w](z)=z^{-n}\Big(-\kappa_{n-1}^{2}C[\Phi_{n-1}^{*}w](z)+1\Big).
\end{equation}
By identifying the 12 and 22 entries in two sides of (4.10), (4.13) and (4.14) respectively follow  from (4.15) and (4.16).
\end{proof}

\begin{rem}
The identities (4.11) and (4.12) are just the classical Szeg\"o recursions. They are equivalent to each other.
\end{rem}

\subsection{The case of OTP} The RH characterization for OTP is uniquely solvable as follows

\begin{thm}[\!\!\cite{dd08}] The RHP (1.3) has a unique solution given by
\begin{equation}
\mathfrak{Y}(z)=\left(
  \begin{array}{cc}
    z^{n}L(\sigma_{n}, \pi_{n})(z) & C[\tau^{-n}L(\sigma_{n}, \pi_{n})w](z) \\
    z^{n}\mathcal{L}(\sigma_{n-1}, \pi_{n-1})(z) & C[\tau^{-n}\mathcal{L}(\sigma_{n-1}, \pi_{n-1})w](z) \\
  \end{array}
\right),
\end{equation}
where \begin{equation}
L(\sigma_{n}, \pi_{n})(z)=\lambda_{1,n}a_{n}\sigma_{n}(z)+\lambda_{2,n}b_{n}\pi_{n}(z),
\end{equation}
\begin{equation}
\mathcal{L}(\sigma_{n-1}, \pi_{n-1})(z)=\lambda_{3,n-1}a_{n-1}\sigma_{n-1}(z)+\lambda_{4,n-1}b_{n-1}\pi_{n-1}(z)
\end{equation}
in which $\sigma_{n}$ and $\pi_{n}$ are the orthonormal Laurent polynomials of the first class on the unit circle with respect to the weight $w$, $a_{n}, b_{n}, \beta_{n}$ are given in (2.11)-(2.13), $\lambda_{1,n}=1$, $\lambda_{2,n}=\beta_{n}+i$, $\lambda_{3,n-1}=-\frac{1}{2}a_{n-1}^{-2}(1+\beta_{n-1} i)$, $\lambda_{4,n-1}=\frac{1}{2}b_{n-1}^{-2}i$, and $C$ is the Cauchy integral operator.
\end{thm}

Set
\begin{equation}
\mathfrak{Y}_{1}(z)=\overline{\mathfrak{Y}\left(\frac{1}{\overline{z}}\right)},\,\,\,\,z\in \mathbb{C}\setminus\partial \mathbb{D},
\end{equation}
then $\mathfrak{Y}_{1}^{+}(t)=\overline{\mathfrak{Y}^{-}(t)}$  and $\mathfrak{Y}_{1}^{-}(t)=\overline{\mathfrak{Y}^{+}(t)}$ for $t\in \partial \mathbb{D}$. Therefore, by the boundary and growth conditions in RHP (1.3), we have
\begin{equation}
\mathfrak{Y}_{1}^{+}(t)=\mathfrak{Y}_{1}^{-}(t)\left(
                   \begin{array}{cc}
                     1 & -t^{2n}w(t) \\
                     0 & 1 \\
                   \end{array}
                 \right),\,\,\, t\in \partial \mathbb{D}
\end{equation}
and
\begin{equation}
\lim_{z\rightarrow 0}\mathfrak{Y}_{1}(z)\left(
                                          \begin{array}{cc}
                                            z^{2n} & 0 \\
                                            0 & z^{-2n+1} \\
                                          \end{array}
                                        \right)
=I.
\end{equation}
Moreover, by straightforward calculations, we have
\begin{equation}
\lim_{z\rightarrow \infty}\mathfrak{Y}_{1}(z)\left(
                                               \begin{array}{cc}
                                                 z & 0 \\
                                                 0 & 1 \\
                                               \end{array}
                                             \right)
=\triangle=\left(
                                                    \begin{array}{cc}
                                                      \triangle_{11} & \triangle_{12} \\
                                                      \triangle_{21} & \triangle_{22} \\
                                                    \end{array}
                                                  \right),
\end{equation}
where
\begin{align}
\triangle_{11}&=-\frac{1}{2}(\imath_{n}+\beta_{n-1}\varsigma_{n}-\zeta_{n})+\frac{i}{2}(\jmath_{n}-\imath_{n}\beta_{n-1}
+\varsigma_{n})\nonumber\\
&=-\alpha_{2n-2}\,\, (\mbox{by Theorem 2.4}),\\
\triangle_{12}&=a_{n}^{2}+b_{n}^{2}(1+\beta_{n}^{2})=\kappa_{2n-1}^{-2}\,\, (\mbox{by Theorem 2.5}),\\
\triangle_{21}&=-\frac{1}{4}\Big(a_{n-1}^{-2}(1+\beta_{n-1}^{2})+b_{n-1}^{-2}\Big)=-\kappa_{2n-2}^{2}\,\, (\mbox{by Theorem 2.3}),\\
\triangle_{22}&=-\frac{1}{2}(\imath_{n}+\varsigma_{n}\beta_{n-1}-\zeta_{n})-\frac{i}{2}(\jmath_{n}-\imath_{n}\beta_{n-1}+\varsigma_{n})=-\overline{\alpha}_{2n-2}.
\end{align}

Let
\begin{align}
\mathfrak{Y}_{2}(z)
=\left(
               \begin{array}{cc}
                     \triangle_{22} & -\triangle_{12} \\
                     \triangle_{21} & -\triangle_{11} \\
                   \end{array}
                 \right)\mathfrak{Y}_{1}(z)\left(
               \begin{array}{cc}
                     z^{2n+1} & 0 \\
                     0 & -z^{-2n+1} \\
                   \end{array}
                 \right).
\end{align}
Noting (or by Theorem 2.7)
\begin{equation}
\det\triangle=|\alpha_{2n-2}|^{2}+\left(\frac{\kappa_{2n-2}}{\kappa_{2n-1}}\right)^{2}=1,
\end{equation}
by (4.23) and (4.28),
\begin{equation}
\lim_{z\rightarrow \infty}\mathfrak{Y}_{2}(z)\left(
               \begin{array}{cc}
                     z^{-2n} & 0 \\
                     0 & z^{2n-1} \\
                   \end{array}
                 \right)=I.
\end{equation}
Thus $\mathfrak{Y}_{2}$ satisfies the following RHP (i.e. the RH characterization (1.3) for OTP)
\begin{equation} (\mbox{RHP for $\mathfrak{Y}_{2}$})\,\,
\begin{cases}
\mathfrak{Y}_{2}\,\, \mbox{is analytic in}\,\,\mathbb{C}\setminus \partial \mathbb{D},\vspace{2mm}\\
\mathfrak{Y}_{2}^{+}(t)=\mathfrak{Y}_{2}^{-}(t)\left(
                   \begin{array}{cc}
                     1 & t^{-2n}w(t) \\
                     0 & 1 \\
                   \end{array}
                 \right)
                 \,\,\mbox{for} \,\,t\in \partial \mathbb{D},\vspace{2mm}\\
\mathfrak{Y}_{2}(z)=\left(I+O(\frac{1}{z})\right)\left(
                                    \begin{array}{cc}
                                      z^{2n} & 0 \\
                                      0 & z^{-2n+1} \\
                                    \end{array}
                                  \right)
                                  \,\,\mbox{as}\,\, z\rightarrow \infty,\vspace{2mm}\\
(\mathfrak{Y}_{2})_{11}(0)=(\mathfrak{Y}_{2})_{21}(0)=0.
\end{cases}
\end{equation}
By the uniqueness of RHP (1.3), $\mathfrak{Y}_{2}(z)=\mathfrak{Y}(z)$ for any $z\in \mathbb{C}\setminus \partial \mathbb{D}$. That is,
\begin{align}
\mathfrak{Y}(z)
=\left(
               \begin{array}{cc}
                     \triangle_{22} & -\triangle_{12} \\
                     \triangle_{21} & -\triangle_{11} \\
                   \end{array}
                 \right)\overline{\mathfrak{Y}\left(\frac{1}{\overline{z}}\right)}\left(
               \begin{array}{cc}
                     z^{2n+1} & 0 \\
                     0 & -z^{-2n+1} \\
                   \end{array}
                 \right),\,\,\,z\in \mathbb{C}\setminus\partial \mathbb{D}.
\end{align}

In order to derive some identities for OTP, we introduce reflectional sets, reflectional and auto-reflectional functions for the unit circle $\partial \mathbb{D}$.

\begin{defn}
A set $\Sigma$ is called a reflectional set for the unit circle $\partial \mathbb{D}$ if both $z\in\Sigma$ and $1/z\in\Sigma$, or
simply a reflectional set, in which $z$ and $1/z$ are called reflection to each other. For example, $\mathbb{C} \setminus\{0\}$ is
a reflectional set for the unit circle.
\end{defn}

\begin{defn}
If $f$ is defined on a reflectional set $\Sigma$, set
\begin{equation}
f_{*}(z)=\overline{f\left(1/\overline{z}\right)},\,\,z\in\Sigma,
\end{equation}
then $f_{*}$ is the reflectional function of $f$ for the unit circle in $\Sigma$, simply reflection.
\end{defn}

\begin{defn}
If $f$ is defined on a reflectional set $\Sigma$ such that
\begin{equation}
f(z)=f_{*}(z),\,\,z\in\Sigma,
\end{equation}
then $f$ is called an auto-reflectional function for the unit circle in $\Sigma$, simply auto-reflection.
\end{defn}

\begin{lem}
Let $\mu$ be a nontrivial probability measure on the unit circle $\partial \mathbb{D}$,
$\{1, \sigma_{n}, \pi_{n}\}$ be the unique system of orthonormal Laurent polynomials of the first class on the unit circle with
respect to $\mu$, then $\sigma_{n}, \pi_{n}$ are auto-reflectional for the unit circle in $\mathbb{C}\setminus\{0\}$.
\end{lem}

\begin{proof}
It immediately follows from that $\displaystyle\frac{z^{n}+z^{-n}}{2}$ and $\displaystyle\frac{z^{n}-z^{-n}}{2i}$ are auto-reflectional for the unit circle in $\mathbb{C}\setminus\{0\}$ and all of the coefficients are real-valued.
\end{proof}

With the above preliminaries, we have
\begin{thm}
Let $\sigma_{n}$, $\pi_{n}$, $a_{n}$, $b_{n}$, $\triangle_{kl}$, $\lambda_{m,n}$ be as above, then
\begin{enumerate}
  \item [($\mathcal{A}$)] The identities
  \begin{align}
  &(\lambda_{1,n}z^{-1}-\overline{\lambda}_{1,n}\triangle_{22})a_{n}\sigma_{n}(z)+(\lambda_{2,n}z^{-1}-\overline{\lambda}_{2,n}\triangle_{22})b_{n}\pi_{n}(z)\nonumber\\
  =&-\triangle_{12}\Big(\overline{\lambda}_{3,n-1}a_{n-1}\sigma_{n-1}(z)+\overline{\lambda}_{4,n-1}b_{n-1}\pi_{n-1}(z)\Big)
  \end{align}
  and
  \begin{align}
  &(\lambda_{3,n-1}z^{-1}+\overline{\lambda}_{3,n-1}\triangle_{11})a_{n-1}\sigma_{n-1}(z)+(\lambda_{4,n-1}z^{-1}+\overline{\lambda}_{4,n-1}\triangle_{11})b_{n-1}\pi_{n-1}(z)\nonumber\\
  =&\triangle_{21}\Big(\overline{\lambda}_{1,n}a_{n}\sigma_{n}(z)+\overline{\lambda}_{2,n}b_{n}\pi_{n}(z)\Big)
  \end{align}
  hold for $z\in \mathbb{C}\setminus\{0\}$;
  \item [($\mathcal{B}$)] The identities
  \begin{align}
  &C[\tau^{-n}(\lambda_{1,n}a_{n}\sigma_{n}+\lambda_{2,n}b_{n}\pi_{n})w](z)\nonumber\\
  =&z^{-2(n-1)}\Big[\triangle_{22}C[\tau^{n-1}(\overline{\lambda}_{1,n}a_{n}\sigma_{n}+\overline{\lambda}_{2,n}b_{n}\pi_{n})w](z)\nonumber\\
  &-\triangle_{12}C[\tau^{n-1}(\overline{\lambda}_{3,n-1}a_{n-1}\sigma_{n-1}+\overline{\lambda}_{4,n-1}b_{n-1}\pi_{n-1})w](z)\Big]
  \end{align}
  and
 \begin{align}
  &C[\tau^{-n}(\lambda_{3,n-1}a_{n-1}\sigma_{n-1}+\lambda_{4,n-1}b_{n-1}\pi_{n-1})w](z)\nonumber\\
  =&z^{-2(n-1)}\Big[\triangle_{21}C[\tau^{n-1}(\overline{\lambda}_{1,n}a_{n}\sigma_{n}+\overline{\lambda}_{2,n}b_{n}\pi_{n})w](z)\nonumber\\
  &-\triangle_{11}C[\tau^{n-1}(\overline{\lambda}_{3,n}a_{n-1}\sigma_{n-1}+\overline{\lambda}_{4,n}b_{n-1}\pi_{n-1})w](z)\Big]
  \end{align}
  hold for $z\in \mathbb{C}\setminus(\partial\mathbb{D}\cup\{0\})$.
  \end{enumerate}
  \end{thm}

  \begin{proof}
  By Lemma 4.8, (4.35) and (4.36) are obtained by identifying the 11 and 21 entries in LHS with the ones in RHS of (4.32). Since
   \begin{equation}
   \overline{C[\tau^{-n}L(\sigma_{n}, \pi_{n})w]\left(\frac{1}{z}\right)}=-zC[\tau^{n-1}\overline{L(\sigma_{n}, \pi_{n})}w](z)
   \end{equation}
   and
   \begin{equation}
   \overline{C[\tau^{-n}\mathcal{L}(\sigma_{n-1}, \pi_{n-1})w]\left(\frac{1}{z}\right)}=-zC[\tau^{n-1}\overline{\mathcal{L}(\sigma_{n-1}, \pi_{n-1})}w](z)
   \end{equation}
   for $z\in \mathbb{C}\setminus(\partial\mathbb{D}\cup\{0\})$,
   then (4.37) and (4.38) follow from comparing the 12 and 22 entries with each other in both sides of (4.32).
  \end{proof}

 In the above theorem, as $z$ is restricted on $\partial \mathbb{D}$, we indeed obtain the following four-term recurrences for OTP.
  \begin{thm}
Let $\sigma_{n}$, $\pi_{n}$, $a_{n}$, $b_{n}$, $\triangle_{kl}$, $\lambda_{m,n}$ be as above, then
\begin{enumerate}
\item [($\mathfrak{A}$)] The identities
\begin{align}
  &(\lambda_{1,n}e^{-i\theta}-\overline{\lambda}_{1,n}\triangle_{22})a_{n}\sigma_{n}(\theta)+(\lambda_{2,n}e^{-i\theta}-\overline{\lambda}_{2,n}\triangle_{22})b_{n}\pi_{n}(\theta)\nonumber\\
  =&-\triangle_{12}\Big(\overline{\lambda}_{3,n-1}a_{n-1}\sigma_{n-1}(\theta)+\overline{\lambda}_{4,n-1}b_{n-1}\pi_{n-1}(\theta)\Big)
  \end{align}
  and
  \begin{align}
  &(\lambda_{3,n-1}e^{-i\theta}+\overline{\lambda}_{3,n-1}\triangle_{11})a_{n-1}\sigma_{n-1}(\theta)+(\lambda_{4,n-1}e^{-i\theta}+\overline{\lambda}_{4,n-1}\triangle_{11})b_{n-1}\pi_{n-1}(\theta)\nonumber\\
  =&\triangle_{21}\Big(\overline{\lambda}_{1,n}a_{n}\sigma_{n}(\theta)+\overline{\lambda}_{2,n}b_{n}\pi_{n}(\theta)\Big)
  \end{align}
  hold for $\theta\in [0, 2\pi)$;
  \item [($\mathfrak{B}$)] The identities
  \begin{align}
  &H[\tau^{-n}(\lambda_{1,n}a_{n}\sigma_{n}+\lambda_{2,n}b_{n}\pi_{n})w](e^{i\theta})\nonumber\\
  =&e^{-i[2(n-1)\theta]}\Big[\triangle_{22}H[\tau^{n-1}(\overline{\lambda}_{1,n}a_{n}\sigma_{n}+\overline{\lambda}_{2,n}b_{n}\pi_{n})w](e^{i\theta})\nonumber\\
  &-\triangle_{12}H[\tau^{n-1}(\overline{\lambda}_{3,n-1}a_{n-1}\sigma_{n-1}+\overline{\lambda}_{4,n-1}b_{n-1}\pi_{n-1})w](e^{i\theta})\Big]
  \end{align}
  and
 \begin{align}
  &H[\tau^{-n}(\lambda_{3,n-1}a_{n-1}\sigma_{n-1}+\lambda_{4,n-1}b_{n-1}\pi_{n-1})w](e^{i\theta})\nonumber\\
  =&e^{-i[2(n-1)\theta]}\Big[\triangle_{21}H[\tau^{n-1}(\overline{\lambda}_{1,n}a_{n}\sigma_{n}+\overline{\lambda}_{2,n}b_{n}\pi_{n})w](e^{i\theta})\nonumber\\
  &-\triangle_{11}H[\tau^{n-1}(\overline{\lambda}_{3,n}a_{n-1}\sigma_{n-1}+\overline{\lambda}_{4,n}b_{n-1}\pi_{n-1})w](e^{i\theta})\Big]
  \end{align}
  hold for $\theta\in [0, 2\pi)$, where $H$ is the Hilbert transform on the unit circle, i.e.
  \begin{equation}
  Hf(t)=P.V. \frac{1}{\pi}\int_{\partial \mathbb{D}}\frac{f(\tau)}{t-\tau}d\tau,\,\,t\in\partial \mathbb{D}
  \end{equation}
  in which $f\in H(\partial \mathbb{D})$.
\end{enumerate}
\end{thm}

\begin{proof}
(4.41) and (4.42) are obvious by identifying $e^{i\theta}\in \partial \mathbb{D}$ with $\theta\in[0,2\pi)$ in (4.35) and (4.36). By the well-known Plemelj formula, viz.
\begin{equation}
C^{\pm}f(t)=\pm\frac{1}{2}f(t)+\frac{i}{2}Hf(t),\,\,t\in \partial \mathbb{D},
\end{equation}
where $f\in H(\partial \mathbb{D})$, taking $z\in \mathbb{D}\rightarrow t=e^{i\theta}$ (or $z\in \mathbb{C}\setminus\overline{\mathbb{D}}\rightarrow t=e^{i\theta}$), then (4.43) and (4.44) easily follow from (4.37), (4.38), (4.41) and (4.42).
\end{proof}

  \begin{rem}
  The identities (4.35) (4.36), (4.41) and (4.42) are four-term recurrences for OTP (exactly, the former two are for OLP of the first class). They are equivalent to each other and also to the ones in \cite{dd08}.
  \end{rem}

  \begin{rem}
 By taking a similar strategy, we can also apply the RH characterization (1.2) for OPUC to derive the above identities for OTP (or OLP of the first class) in Theorems 4.9 and 4.10 when $2n-1$ is in place of $n$ as stated in the Introduction.
  \end{rem}

  \begin{rem}
  By the mutual representation theorem (Theorem 2.1) and Theorem 4.2, we can directly get some identities for OTP in different forms. They are equivalent to (4.35)-(4.38) and (4.41)-(4.44). By this approach, the four-term recurrences (corresponding to (4.35), (4.36), (4.41) and (4.42)) were obtained in \cite{dd08}.
  \end{rem}

\section{A Strong Favard theorem}

Theorem 3.1 tells us that there exist many nontrivial probability measures $d\mu$ corresponding to any fixed system of three-tuples $\{(a_{n}^{(0)},b_{n}^{(0)},\beta_{n}^{(0)})\}$ satisfying (3.1). That is to say that the system of three-tuples $\{(a_{n}^{(0)},b_{n}^{(0)},\beta_{n}^{(0)})\}$ with (3.1) is not sufficient to uniquely determine the nontrivial probability measure $d\mu$.
As stated in Remark 3.2, we need consider a system of seven-tuples $\{(a_{n}^{(0)}, b_{n}^{(0)},\beta_{n}^{(0)},$ $\imath_{n}^{(0)},\jmath_{n}^{(0)},\varsigma_{n}^{(0)},\zeta_{n}^{(0)})\}$ with some suitable properties in order to uniquely determine the nontrivial probability measure $d\mu$. In what follows, we will discuss it in detail.

At first, we give some more relations on the coefficients $a_{n}, b_{n},\beta_{n},\imath_{n},\jmath_{n},\varsigma_{n},\zeta_{n}$ of OTP and $\alpha_{n}, \kappa_{n}$ of OPUC. To this end, the following basic facts are required.
\begin{lem}
Let $\Phi_{n}$ be the monic orthogonal polynomial on the unit circle of order $n$ with respect to $\mu$, and $\Phi^{*}_{n}$ be the reversed polynomial of $\Phi_{n}$, then
\begin{equation}
\langle 1,z\Phi_{n}^{*}\rangle_{\mathbb{R}}=\int_{\partial \mathbb{D}}z\Phi_{n}^{*}(\tau)d\mu(\tau)=-a_{n+1,n}
\end{equation}
and
\begin{equation}
\langle 1,z^{2}\Phi_{n}\rangle_{\mathbb{R}}=\int_{\partial \mathbb{D}}\tau^{2}\Phi_{n}(\tau)d\mu(\tau)=\alpha_{n}^{-1}\Big(a_{n+2,n+1}-a_{n+1,n}\Big),
\end{equation}
where the Verblunsky coefficient $\alpha_{n}$ is restricted in $\mathbb{D}\setminus\{0\}$, and $a_{n+1,n}$ is given by
\begin{equation}
\Phi_{n+1}(z)=z^{n+1}+a_{n+1,n}z^{n}+\mbox{lower order}.
\end{equation}
\end{lem}
\begin{proof}
By (2.6), we have
\begin{align}
\int_{\partial \mathbb{D}}z\Phi_{n}^{*}(\tau)d\mu(\tau)=\int_{\partial \mathbb{D}}\tau^{n+1}\overline{\Phi_{n}(\tau)}d\mu(\tau)=\langle z^{n+1},\Phi_{n}\rangle_{\mathbb{C}}.
\end{align}
Thus (5.1) immediately follows from (5.3) by the orthogonality. For $\alpha_{n}\in\mathbb{D}\setminus\{0\}$, by Szeg\"o recurrence,
\begin{align}
\int_{\partial \mathbb{D}}\tau^{2}\Phi_{n}(\tau)d\mu(\tau)=\alpha_{n}^{-1}\Big[\int_{\partial \mathbb{D}}\tau\Phi_{n}^{*}(\tau)d\mu(\tau)-\int_{\partial \mathbb{D}}\tau\Phi_{n+1}^{*}(\tau)d\mu(\tau)\Big].
\end{align}
So (5.2) holds by applying (5.1).
\end{proof}

\begin{thm} Let $a_{n}, b_{n},\beta_{n},\imath_{n},\jmath_{n},\varsigma_{n},\zeta_{n},\alpha_{n},\kappa_{n}$ be given in Section 2, then
\begin{align}
\imath_{n+1}=&\frac{1}{4}\Lambda_{n}^{-1}a_{n}^{-2}b_{n}^{-2}\Big(\alpha_{2n-1}^{-1}+1\Big)a_{2n,2n-1}i-\frac{1}{4}\Big[\Lambda_{n}^{-1}a_{n}^{-2}b_{n}^{-2}\alpha_{2n-1}^{-1}i\nonumber\\
&+a_{n}^{-2}(1-\beta_{n}i)\Big]a_{2n+1,2n}
+\frac{1}{4}a_{n}^{-2}(1-\beta_{n}i)\overline{\alpha}_{2n}^{-1}\Big(\kappa_{2n}^{-2}-\kappa_{2n+1}^{-2}\Big),
\end{align}
\begin{align}
\jmath_{n+1}=&\frac{1}{4}\Lambda_{n}^{-1}a_{n}^{-2}b_{n}^{-2}(1+\beta_{n}i)\Big(\alpha_{2n-1}^{-1}+1\Big)a_{2n,2n-1}-\frac{1}{4}\Big[\Lambda_{n}^{-1}a_{n}^{-2}b_{n}^{-2}(1+\beta_{n}i)\alpha_{2n-1}^{-1}\nonumber\\
&+b_{n}^{-2}i\Big]a_{2n+1,2n}
+\frac{1}{4}b_{n}^{-2}\overline{\alpha}_{2n}^{-1}\Big(\kappa_{2n}^{-2}-\kappa_{2n+1}^{-2}\Big)i,
\end{align}
\begin{align}
\varsigma_{n+1}=&\frac{1}{4}\Lambda_{n}^{-1}a_{n}^{-2}b_{n}^{-2}\Big(\alpha_{2n-1}^{-1}-1\Big)a_{2n,2n-1}-\frac{1}{4i}\Big[\Lambda_{n}^{-1}a_{n}^{-2}b_{n}^{-2}\alpha_{2n-1}^{-1}i\nonumber\\
&+a_{n}^{-2}(1-\beta_{n}i)\Big]a_{2n+1,2n}
-\frac{1}{4i}a_{n}^{-2}(1-\beta_{n}i)\overline{\alpha}_{2n}^{-1}\Big(\kappa_{2n}^{-2}-\kappa_{2n+1}^{-2}\Big)
\end{align}
and
\begin{align}
\zeta_{n+1}=&\frac{1}{4i}\Lambda_{n}^{-1}a_{n}^{-2}b_{n}^{-2}(1+\beta_{n}i)\Big(\alpha_{2n-1}^{-1}-1\Big)a_{2n,2n-1}-\frac{1}{4i}\Big[\Lambda_{n}^{-1}a_{n}^{-2}b_{n}^{-2}(1+\beta_{n}i)\alpha_{2n-1}^{-1}\nonumber\\
&+b_{n}^{-2}i\Big]a_{2n+1,2n}
-\frac{1}{4}b_{n}^{-2}\overline{\alpha}_{2n}^{-1}\Big(\kappa_{2n}^{-2}-\kappa_{2n+1}^{-2}\Big)
\end{align}
for $n\in \mathbb{N}\cup\{0\}$, where $\alpha_{2n-1},\alpha_{2n}\in \mathbb{D}\setminus\{0\}$, $a_{2n+1,2n}, a_{2n,2n-1}$ are given by (5.3).
\end{thm}
\begin{proof}
It is enough to prove (5.6) and (5.7). Similar are to (5.8) and (5.9).
By Theorem 2.2, Lemmas 5.1 and 2.8 and Szeg\"o recurrence, when $\alpha_{2n-1},\alpha_{2n}\in \mathbb{D}\setminus\{0\}$, we have
\begin{align}
\langle z^{n+1},a_{n}\sigma_{n}\rangle_{\mathbb{R}}=&\langle z^{n+1},-\frac{1}{2}z^{-n}[\Lambda_{n}^{-1}b_{n}^{-2}iz\Phi_{2n-1}-(1-\beta_{n}i)\Phi_{2n}^{*}]\rangle_{\mathbb{R}}\nonumber\\
=&-\frac{1}{2}\Lambda_{n}^{-1}b_{n}^{-2}i\langle 1,z^{2}\Phi_{2n-1}\rangle_{\mathbb{R}}+\frac{1}{2}(1-\beta_{n}i)\langle 1,z\Phi_{2n}^{*}\rangle_{\mathbb{R}}\nonumber\\
=&\frac{1}{2}\Lambda_{n}^{-1}b_{n}^{-2}\alpha_{2n-1}^{-1}a_{2n,2n-1}i-\frac{1}{2}\Big[\Lambda_{n}^{-1}b_{n}^{-2}\alpha_{2n-1}^{-1}i\nonumber\\
&+(1-\beta_{n}i)\Big]a_{2n+1,2n},
\end{align}
\begin{align}
\langle z^{-(n+1)},a_{n}\sigma_{n}\rangle_{\mathbb{R}}=&\langle z^{-(n+1)},-\frac{1}{2}z^{-n}[\Lambda_{n}^{-1}b_{n}^{-2}iz\Phi_{2n-1}-(1-\beta_{n}i)\Phi_{2n}^{*}]\rangle_{\mathbb{R}}\nonumber\\
=&-\frac{1}{2}\Lambda_{n}^{-1}b_{n}^{-2}i\langle 1,z^{-2n}\Phi_{2n-1}\rangle_{\mathbb{R}}+\frac{1}{2}(1-\beta_{n}i)\langle 1,z^{-(2n+1)}\Phi_{2n}^{*}\rangle_{\mathbb{R}}\nonumber\\
=&-\frac{1}{2}\Lambda_{n}^{-1}b_{n}^{-2}i\langle z^{2n},\Phi_{2n-1}\rangle_{\mathbb{C}}+\frac{1}{2}(1-\beta_{n}i)\overline{\langle 1,z\Phi_{2n}\rangle}_{\mathbb{R}}\nonumber\\
=&\frac{1}{2}\Lambda_{n}^{-1}b_{n}^{-2}a_{2n,2n-1}i+\frac{1}{2}(1-\beta_{n}i)\overline{\alpha}_{2n}^{-1}\Big(\kappa_{2n}^{-2}-\kappa_{2n+1}^{-2}\Big),
\end{align}
\begin{align}
\langle z^{n+1},b_{n}\pi_{n}\rangle_{\mathbb{R}}=&\langle z^{n+1},-\frac{1}{2}z^{-n}[\Lambda_{n}^{-1}a_{n}^{-2}(1+\beta_{n}i)z\Phi_{2n-1}(z)-i\Phi_{2n}^{*}(z)]\rangle_{\mathbb{R}}\nonumber\\
=&-\frac{1}{2}\Lambda_{n}^{-1}a_{n}^{-2}(1+\beta_{n}i)\langle 1,z^{2}\Phi_{2n-1}\rangle_{\mathbb{R}}+\frac{i}{2}\langle 1,z\Phi_{2n}^{*}\rangle_{\mathbb{R}}\nonumber\\
=&\frac{1}{2}\Lambda_{n}^{-1}a_{n}^{-2}(1+\beta_{n}i)\alpha_{2n-1}^{-1}a_{2n,2n-1}-\frac{1}{2}\Big[\Lambda_{n}^{-1}a_{n}^{-2}(1+\beta_{n}i)\alpha_{2n-1}^{-1}\nonumber\\&+i\Big]a_{2n+1,2n}
\end{align}
and
\begin{align}
\langle z^{-(n+1)},b_{n}\pi_{n}\rangle_{\mathbb{R}}=&\langle z^{-(n+1)},-\frac{1}{2}z^{-n}[\Lambda_{n}^{-1}a_{n}^{-2}(1+\beta_{n}i)z\Phi_{2n-1}(z)-i\Phi_{2n}^{*}(z)]\rangle_{\mathbb{R}}\nonumber\\
=&-\frac{1}{2}\Lambda_{n}^{-1}a_{n}^{-2}(1+\beta_{n}i)\langle 1,z^{-2n}\Phi_{2n-1}\rangle_{\mathbb{R}}+\frac{i}{2}\langle 1,z^{-(2n+1)}\Phi_{2n}^{*}\rangle_{\mathbb{R}}\nonumber\\
=&-\frac{1}{2}\Lambda_{n}^{-1}a_{n}^{-2}(1+\beta_{n}i)\langle z^{2n},\Phi_{2n-1}\rangle_{\mathbb{C}}+\frac{i}{2}\overline{\langle 1,z\Phi_{2n}\rangle}_{\mathbb{R}}\nonumber\\
=&\frac{1}{2}\Lambda_{n}^{-1}a_{n}^{-2}(1+\beta_{n}i)a_{2n,2n-1}+\frac{i}{2}\overline{\alpha}_{2n}^{-1}\Big(\kappa_{2n}^{-2}-\kappa_{2n+1}^{-2}\Big).
\end{align}

Thus
\begin{align}
\imath_{n+1}=&\langle\frac{z^{n+1}+z^{-(n+1)}}{2},a_{n}^{-1}\sigma_{n}\rangle_{\mathbb{R}}=a_{n}^{-2}\langle\frac{z^{n+1}+z^{-(n+1)}}{2},a_{n}\sigma_{n}\rangle_{\mathbb{R}}\nonumber\\
=&\frac{1}{4}\Lambda_{n}^{-1}a_{n}^{-2}b_{n}^{-2}\Big(\alpha_{2n-1}^{-1}+1\Big)a_{2n,2n-1}i-\frac{1}{4}\Big[\Lambda_{n}^{-1}a_{n}^{-2}b_{n}^{-2}\alpha_{2n-1}^{-1}i\nonumber\\
&+a_{n}^{-2}(1-\beta_{n}i)\Big]a_{2n+1,2n}
+\frac{1}{4}a_{n}^{-2}(1-\beta_{n}i)\overline{\alpha}_{2n}^{-1}\Big(\kappa_{2n}^{-2}-\kappa_{2n+1}^{-2}\Big)
\end{align}
and
\begin{align}
\jmath_{n+1}=&\langle\frac{z^{n+1}+z^{-(n+1)}}{2},b_{n}^{-1}\pi_{n}\rangle_{\mathbb{R}}=b_{n}^{-2}\langle\frac{z^{n+1}+z^{-(n+1)}}{2},b_{n}\pi_{n}\rangle_{\mathbb{R}}\nonumber\\
=&\frac{1}{4}\Lambda_{n}^{-1}a_{n}^{-2}b_{n}^{-2}(1+\beta_{n}i)\Big(\alpha_{2n-1}^{-1}+1\Big)a_{2n,2n-1}-\frac{1}{4}\Big[\Lambda_{n}^{-1}a_{n}^{-2}b_{n}^{-2}(1+\beta_{n}i)\alpha_{2n-1}^{-1}\nonumber\\
&+b_{n}^{-2}i\Big]a_{2n+1,2n}
+\frac{1}{4}b_{n}^{-2}\overline{\alpha}_{2n}^{-1}\Big(\kappa_{2n}^{-2}-\kappa_{2n+1}^{-2}\Big)i,
\end{align}
where $\alpha_{2n-1},\alpha_{2n}\in \mathbb{D}\setminus\{0\}$.
\end{proof}

Denote
\begin{equation}
A=\left(
    \begin{array}{cc}
      \frac{1}{4}\Lambda_{n}^{-1}a_{n}^{-2}b_{n}^{-2}\Big(\alpha_{2n-1}^{-1}+1\Big)i & -\frac{1}{4}\Big[\Lambda_{n}^{-1}a_{n}^{-2}b_{n}^{-2}\alpha_{2n-1}^{-1}i+a_{n}^{-2}(1-\beta_{n}i)\Big] \vspace{2mm}\\
      \frac{1}{4}\Lambda_{n}^{-1}a_{n}^{-2}b_{n}^{-2}(1+\beta_{n}i)\Big(\alpha_{2n-1}^{-1}+1\Big) & -\frac{1}{4}\Big[\Lambda_{n}^{-1}a_{n}^{-2}b_{n}^{-2}(1+\beta_{n}i)\alpha_{2n-1}^{-1}+b_{n}^{-2}i\Big] \\
    \end{array}
  \right),
\end{equation}
then by (2.21) and (2.37),
\begin{align}
|A|=&\left|\begin{array}{cc}
      \frac{1}{4}\Lambda_{n}^{-1}a_{n}^{-2}b_{n}^{-2}\Big(\alpha_{2n-1}^{-1}+1\Big)i & -\frac{1}{4}\Lambda_{n}^{-1}a_{n}^{-2}b_{n}^{-2}\alpha_{2n-1}^{-1}i\vspace{2mm}\\
      \frac{1}{4}\Lambda_{n}^{-1}a_{n}^{-2}b_{n}^{-2}(1+\beta_{n}i)\Big(\alpha_{2n-1}^{-1}+1\Big) & -\frac{1}{4}\Lambda_{n}^{-1}a_{n}^{-2}b_{n}^{-2}(1+\beta_{n}i)\alpha_{2n-1}^{-1} \\
    \end{array}
\right|\nonumber\\
&+\left|\begin{array}{cc}
      \frac{1}{4}\Lambda_{n}^{-1}a_{n}^{-2}b_{n}^{-2}\Big(\alpha_{2n-1}^{-1}+1\Big)i & -\frac{1}{4}a_{n}^{-2}(1-\beta_{n}i)\vspace{2mm}\\
      \frac{1}{4}\Lambda_{n}^{-1}a_{n}^{-2}b_{n}^{-2}(1+\beta_{n}i)\Big(\alpha_{2n-1}^{-1}+1\Big) & -\frac{1}{4}b_{n}^{-2}i \\
    \end{array}
\right|\nonumber\\
=&-\frac{1}{16}\Lambda_{n}^{-1}a_{n}^{-2}b_{n}^{-2}\Big(\alpha_{2n-1}^{-1}+1\Big)\left|\begin{array}{cc}
      i & a_{n}^{-2}(1-\beta_{n}i)\vspace{2mm}\\
      (1+\beta_{n}i) & b_{n}^{-2}i \\
    \end{array}
\right|\nonumber\\
=&\frac{1}{16}\Lambda_{n}^{-1}a_{n}^{-2}b_{n}^{-2}\Big(\alpha_{2n-1}^{-1}+1\Big)[a_{n}^{-2}(1+\beta_{n}^{2})+b_{n}^{-2}]\nonumber\\
=&\frac{1}{8}a_{n}^{-2}b_{n}^{-2}\Big(\alpha_{2n-1}^{-1}+1\Big)i\neq 0.
\end{align}
So $A$ is invertible.

By Theorem 5.2, we can represent the coefficient $a_{n+1,n}$ of OPUC by the coefficients $a_{n}, b_{n},\beta_{n},\imath_{n},\jmath_{n},\varsigma_{n},\zeta_{n}$ of OTP as follows.
\begin{thm}
Let $a_{n}, b_{n},\beta_{n},\imath_{n},\jmath_{n},\varsigma_{n},\zeta_{n},\alpha_{n},\kappa_{n},a_{n+1,n}$ be as above, then
\begin{equation}
\left(
  \begin{array}{c}
    a_{2n,2n-1} \vspace{1mm}\\
    a_{2n+1,2n} \\
  \end{array}
\right)=A^{-1}\left(
                \begin{array}{c}
                  \imath_{n+1}-\frac{1}{4}a_{n}^{-2}(1-\beta_{n}i)\overline{\alpha}_{2n}^{-1}\Big(\kappa_{2n}^{-2}-\kappa_{2n+1}^{-2}\Big) \vspace{1mm}\\
                  \jmath_{n+1}-\frac{1}{4}b_{n}^{-2}\overline{\alpha}_{2n}^{-1}\Big(\kappa_{2n}^{-2}-\kappa_{2n+1}^{-2}\Big)i \\
                \end{array}
              \right),
\end{equation}
where $A$ is given by (5.16).
\end{thm}
\begin{proof}
It follows from (5.6), (5.7) and the invertibility of $A$.
\end{proof}
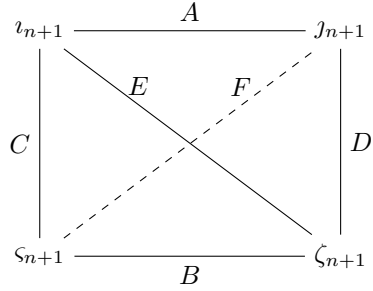
\begin{center}
\begin{tikzpicture}
\node (a) at (0,0) {$\varsigma_{n+1}$};
\node (b) at (0,3) {$\imath_{n+1}$};
\node (c) at (4,0) {$\zeta_{n+1}$};
\node (d) at (4,3) {$\jmath_{n+1}$};
\draw (a) -- node[left]{$C$} (b);
\draw (a) -- node[below]{$B$} (c);
\draw (b) -- node[above]{$A$} (d);
\draw (c) -- node[right]{$D$} (d);
\draw [dashed] (a) -- node[above,pos=0.7]{$F$} (d);
\draw (b) -- node[above,pos=0.3]{$E$} (c);
\end{tikzpicture}
\end{center}
\begin{center}
Derivation of $a_{n+1,n}$ from different ways
\end{center}

\begin{rem}
In Theorem 5.3, we derive $a_{2n+1,2n}$ and $a_{2n,2n-1}$ in terms of $\imath_{n+1}$ and $\jmath_{n+1}$. In fact, there are many different ways shown in the above figure from which to deduce them. For instance, like as $A$, let $B$ be the coefficient matrix for $a_{2n+1,2n}$ and $a_{2n,2n-1}$ in (5.8) and (5.9). Namely,
\begin{equation}
B=\left(
    \begin{array}{cc}
      \frac{1}{4}\Lambda_{n}^{-1}a_{n}^{-2}b_{n}^{-2}\Big(\alpha_{2n-1}^{-1}-1\Big) & -\frac{1}{4i}\Big[\Lambda_{n}^{-1}a_{n}^{-2}b_{n}^{-2}\alpha_{2n-1}^{-1}i+a_{n}^{-2}(1-\beta_{n}i)\Big] \vspace{2mm}\\
      \frac{1}{4i}\Lambda_{n}^{-1}a_{n}^{-2}b_{n}^{-2}(1+\beta_{n}i)\Big(\alpha_{2n-1}^{-1}-1\Big) & -\frac{1}{4i}\Big[\Lambda_{n}^{-1}a_{n}^{-2}b_{n}^{-2}(1+\beta_{n}i)\alpha_{2n-1}^{-1}+b_{n}^{-2}i\Big] \\
    \end{array}
  \right).
\end{equation}
Then we can use $B$ to express $a_{2n+1,2n}$ and $a_{2n,2n-1}$ in terms of $\varsigma_{n+1}$ and $\zeta_{n+1}$. So are to use $C,D$ and $E$ in the solid-line cases shown in the above figure. However, a further condition will be required for $F$ in the dashed-line case. Here $C, D, E, F$ have the similar sense as $A$ and $B$. Such observations are based on the following results about the evaluations for the determinants of these coefficient matrices.
\end{rem}

\begin{thm} Let $B, C, D, E, F$ be the coefficient matrices stated in the above remark, then
\begin{equation}
|B|=-\frac{1}{8}a_{n}^{-2}b_{n}^{-2}\Big(\alpha_{2n-1}^{-1}-1\Big)i,
\end{equation}
\begin{equation}
|C|=\frac{1}{8}\Lambda_{n}^{-1}a_{n}^{-4}b_{n}^{-2}\alpha_{2n-1}^{-1}(1+\beta_{n}i),
\end{equation}
\begin{equation}
|D|=-\frac{1}{8}\Lambda_{n}^{-1}a_{n}^{-2}b_{n}^{-4}\alpha_{2n-1}^{-1},
\end{equation}
\begin{equation}
|E|=\frac{1}{8i}\Lambda_{n}^{-1}a_{n}^{-2}b_{n}^{-2}\alpha_{2n-1}^{-1}\Big[a_{n}^{-2}+b_{n}^{-2}+a_{n}^{-2}\beta_{n}i\Big]
\end{equation}
and
\begin{equation}
|F|=-\frac{1}{8i}\Lambda_{n}^{-1}a_{n}^{-4}b_{n}^{-2}\alpha_{2n-1}^{-1}\beta_{n}(\beta_{n}-i).
\end{equation}
\end{thm}
\begin{proof}
By applying (2.21) and (2.37) as well as basic properties of determinants, we have
\begin{align}
|B|=&\left|\begin{array}{cc}
      \frac{1}{4}\Lambda_{n}^{-1}a_{n}^{-2}b_{n}^{-2}\Big(\alpha_{2n-1}^{-1}-1\Big) & -\frac{1}{4}\Lambda_{n}^{-1}a_{n}^{-2}b_{n}^{-2}\alpha_{2n-1}^{-1}\vspace{2mm}\\
      \frac{1}{4i}\Lambda_{n}^{-1}a_{n}^{-2}b_{n}^{-2}(1+\beta_{n}i)\Big(\alpha_{2n-1}^{-1}-1\Big) & -\frac{1}{4i}\Lambda_{n}^{-1}a_{n}^{-2}b_{n}^{-2}(1+\beta_{n}i)\alpha_{2n-1}^{-1} \\
    \end{array}
\right|\nonumber\\
&+\left|\begin{array}{cc}
      \frac{1}{4}\Lambda_{n}^{-1}a_{n}^{-2}b_{n}^{-2}\Big(\alpha_{2n-1}^{-1}-1\Big) & -\frac{1}{4i}a_{n}^{-2}(1-\beta_{n}i)\vspace{2mm}\\
      \frac{1}{4i}\Lambda_{n}^{-1}a_{n}^{-2}b_{n}^{-2}(1+\beta_{n}i)\Big(\alpha_{2n-1}^{-1}-1\Big) & -\frac{1}{4}b_{n}^{-2}\\
    \end{array}
\right|\nonumber\\
=&\frac{1}{16}\Lambda_{n}^{-1}a_{n}^{-2}b_{n}^{-2}\Big(\alpha_{2n-1}^{-1}-1\Big)\left|\begin{array}{cc}
      i & a_{n}^{-2}(1-\beta_{n}i)\vspace{2mm}\\
      (1+\beta_{n}i) & b_{n}^{-2}i \\
    \end{array}
\right|\nonumber\\
=&-\frac{1}{16}\Lambda_{n}^{-1}a_{n}^{-2}b_{n}^{-2}\Big(\alpha_{2n-1}^{-1}-1\Big)[a_{n}^{-2}(1+\beta_{n}^{2})+b_{n}^{-2}]\nonumber\\
=&-\frac{1}{8}a_{n}^{-2}b_{n}^{-2}\Big(\alpha_{2n-1}^{-1}-1\Big)i,\nonumber
\end{align}
\begin{align}
|C|=&\left|\begin{array}{cc}
      \frac{1}{4}\Lambda_{n}^{-1}a_{n}^{-2}b_{n}^{-2}\Big(\alpha_{2n-1}^{-1}+1\Big)i & -\frac{1}{4}\Big[\Lambda_{n}^{-1}a_{n}^{-2}b_{n}^{-2}\alpha_{2n-1}^{-1}i+a_{n}^{-2}(1-\beta_{n}i)\Big]\vspace{2mm}\\
      \frac{1}{4}\Lambda_{n}^{-1}a_{n}^{-2}b_{n}^{-2}\Big(\alpha_{2n-1}^{-1}-1\Big) & -\frac{1}{4i}\Big[\Lambda_{n}^{-1}a_{n}^{-2}b_{n}^{-2}\alpha_{2n-1}^{-1}i+a_{n}^{-2}(1-\beta_{n}i)\Big] \\
    \end{array}
\right|\nonumber\\
=&-\frac{1}{16}\Lambda_{n}^{-1}a_{n}^{-2}b_{n}^{-2}\Big[\Lambda_{n}^{-1}a_{n}^{-2}b_{n}^{-2}\alpha_{2n-1}^{-1}i+a_{n}^{-2}(1-\beta_{n}i)\Big]\left|\begin{array}{cc}
      \alpha_{2n-1}^{-1}+1 & 1\vspace{2mm}\\
      \alpha_{2n-1}^{-1}-1 & 1 \\
    \end{array}
\right|\nonumber\\
=&-\frac{1}{8}\Lambda_{n}^{-1}a_{n}^{-4}b_{n}^{-2}\alpha_{2n-1}^{-1}\Big[\Lambda_{n}^{-1}b_{n}^{-2}i+\alpha_{2n-1}(1-\beta_{n}i)\Big]\nonumber\\
=&\frac{1}{8}\Lambda_{n}^{-1}a_{n}^{-4}b_{n}^{-2}\alpha_{2n-1}^{-1}(1+\beta_{n}i),\nonumber
\end{align}
\begin{align}
|D|=&\left|\begin{array}{cc}
      \frac{1}{4}\Lambda_{n}^{-1}a_{n}^{-2}b_{n}^{-2}(1+\beta_{n}i)\Big(\alpha_{2n-1}^{-1}+1\Big) & -\frac{1}{4}\Big[\Lambda_{n}^{-1}a_{n}^{-2}b_{n}^{-2}(1+\beta_{n}i)\alpha_{2n-1}^{-1}+b_{n}^{-2}i\Big]\vspace{2mm}\\
      \frac{1}{4i}\Lambda_{n}^{-1}a_{n}^{-2}b_{n}^{-2}(1+\beta_{n}i)\Big(\alpha_{2n-1}^{-1}-1\Big) & -\frac{1}{4i}\Big[\Lambda_{n}^{-1}a_{n}^{-2}b_{n}^{-2}(1+\beta_{n}i)\alpha_{2n-1}^{-1}+b_{n}^{-2}i\Big]\\
    \end{array}
\right|\nonumber\\
=&-\frac{1}{16i}\Lambda_{n}^{-1}a_{n}^{-2}b_{n}^{-2}\Big[\Lambda_{n}^{-1}a_{n}^{-2}b_{n}^{-2}(1+\beta_{n}i)\alpha_{2n-1}^{-1}+b_{n}^{-2}i\Big]\left|\begin{array}{cc}
      \alpha_{2n-1}^{-1}+1 & 1\vspace{2mm}\\
      \alpha_{2n-1}^{-1}-1 & 1 \\
    \end{array}
\right|\nonumber\\
=&-\frac{1}{8i}\Lambda_{n}^{-1}a_{n}^{-2}b_{n}^{-4}\alpha_{2n-1}^{-1}\Big[\Lambda_{n}^{-1}a_{n}^{-2}(1+\beta_{n}i)+\alpha_{2n-1}i\Big]\nonumber\\
=&-\frac{1}{8}\Lambda_{n}^{-1}a_{n}^{-2}b_{n}^{-4}\alpha_{2n-1}^{-1},\nonumber
\end{align}
\begin{align}
|E|=&\left|\begin{array}{cc}
      \frac{1}{4}\Lambda_{n}^{-1}a_{n}^{-2}b_{n}^{-2}\Big(\alpha_{2n-1}^{-1}+1\Big)i & -\frac{1}{4}\Lambda_{n}^{-1}a_{n}^{-2}b_{n}^{-2}\alpha_{2n-1}^{-1}i\vspace{2mm}\\
      \frac{1}{4i}\Lambda_{n}^{-1}a_{n}^{-2}b_{n}^{-2}(1+\beta_{n}i)\Big(\alpha_{2n-1}^{-1}-1\Big) & -\frac{1}{4i}\Lambda_{n}^{-1}a_{n}^{-2}b_{n}^{-2}(1+\beta_{n}i)\alpha_{2n-1}^{-1} \\
    \end{array}
\right|\nonumber\\
&+\left|\begin{array}{cc}
      \frac{1}{4}\Lambda_{n}^{-1}a_{n}^{-2}b_{n}^{-2}\Big(\alpha_{2n-1}^{-1}+1\Big)i & -\frac{1}{4}a_{n}^{-2}(1-\beta_{n}i)\vspace{2mm}\\
      \frac{1}{4i}\Lambda_{n}^{-1}a_{n}^{-2}b_{n}^{-2}(1+\beta_{n}i)\Big(\alpha_{2n-1}^{-1}-1\Big) & -\frac{1}{4}b_{n}^{-2} \\
    \end{array}
\right|\nonumber
\end{align}
\begin{align}
=&-\frac{1}{16}(\Lambda_{n}^{-1}a_{n}^{-2}b_{n}^{-2})^{2}(1+\beta_{n}i)\alpha_{2n-1}^{-1}\left|\begin{array}{cc}
      \alpha_{2n-1}^{-1}+1 & 1\vspace{2mm}\\
      \alpha_{2n-1}^{-1}-1 & 1 \\
    \end{array}
\right|\nonumber\\
&-\frac{1}{16i}\Lambda_{n}^{-1}a_{n}^{-2}b_{n}^{-2}\left|\begin{array}{cc}
      (\alpha_{2n-1}^{-1}+1)i & a_{n}^{-2}(1-\beta_{n}i)\vspace{2mm}\\
      (\alpha_{2n-1}^{-1}-1)(1+\beta_{n}i) & b_{n}^{-2}i \\
    \end{array}
\right|\nonumber\\
=&-\frac{1}{8}(\Lambda_{n}^{-1}a_{n}^{-2}b_{n}^{-2})^{2}(1+\beta_{n}i)\alpha_{2n-1}^{-1}+\frac{1}{16i}\Lambda_{n}^{-1}a_{n}^{-2}b_{n}^{-2}\alpha_{2n-1}^{-1}[a_{n}^{-2}(1+\beta_{n}^{2})+b_{n}^{-2}]\nonumber\\
&+\frac{1}{16i}\Lambda_{n}^{-1}a_{n}^{-2}b_{n}^{-2}[b_{n}^{-2}-a_{n}^{-2}(1+\beta_{n}^{2})]\nonumber\\
=&\frac{1}{16i}\Lambda_{n}^{-1}a_{n}^{-2}b_{n}^{-2}\alpha_{2n-1}^{-1}\Big\{\kappa_{2n}^{-2}a_{n}^{-2}b_{n}^{-2}(1+\beta_{n}i)+[a_{n}^{-2}(1+\beta_{n}^{2})+b_{n}^{-2}]\nonumber\\
&+\alpha_{2n-1}[b_{n}^{-2}-a_{n}^{-2}(1+\beta_{n}^{2})]\Big\}\nonumber\\
=&\frac{1}{8i}\Lambda_{n}^{-1}a_{n}^{-2}b_{n}^{-2}\alpha_{2n-1}^{-1}[a_{n}^{-2}(1+\beta_{n}i)+b_{n}^{-2}]\nonumber
\end{align}
and
\begin{align}
|F|=&\left|\begin{array}{cc}
      \frac{1}{4}\Lambda_{n}^{-1}a_{n}^{-2}b_{n}^{-2}(1+\beta_{n}i)\Big(\alpha_{2n-1}^{-1}+1\Big) & -\frac{1}{4}\Lambda_{n}^{-1}a_{n}^{-2}b_{n}^{-2}(1+\beta_{n}i)\alpha_{2n-1}^{-1}\vspace{2mm}\\
      \frac{1}{4}\Lambda_{n}^{-1}a_{n}^{-2}b_{n}^{-2}\Big(\alpha_{2n-1}^{-1}-1\Big) & -\frac{1}{4}\Lambda_{n}^{-1}a_{n}^{-2}b_{n}^{-2}\alpha_{2n-1}^{-1} \\
    \end{array}
\right|\nonumber\\
&+\left|\begin{array}{cc}
      \frac{1}{4}\Lambda_{n}^{-1}a_{n}^{-2}b_{n}^{-2}(1+\beta_{n}i)\Big(\alpha_{2n-1}^{-1}+1\Big) & -\frac{1}{4}b_{n}^{-2}i\vspace{2mm}\\
      \frac{1}{4}\Lambda_{n}^{-1}a_{n}^{-2}b_{n}^{-2}\Big(\alpha_{2n-1}^{-1}-1\Big) & -\frac{1}{4i}a_{n}^{-2}(1-\beta_{n}i)\\
    \end{array}
\right|\nonumber\\
=&-\frac{1}{16}(\Lambda_{n}^{-1}a_{n}^{-2}b_{n}^{-2})^{2}(1+\beta_{n}i)\alpha_{2n-1}^{-1}\left|\begin{array}{cc}
      \alpha_{2n-1}^{-1}+1 & 1\vspace{2mm}\\
      \alpha_{2n-1}^{-1}-1 & 1 \\
    \end{array}
\right|\nonumber\\
&-\frac{1}{16i}\Lambda_{n}^{-1}a_{n}^{-2}b_{n}^{-2}\left|\begin{array}{cc}
      (\alpha_{2n-1}^{-1}+1)(1+\beta_{n}i) & b_{n}^{-2}i\vspace{2mm}\\
      (\alpha_{2n-1}^{-1}-1)i & a_{n}^{-2}(1-\beta_{n}i) \\
    \end{array}
\right|\nonumber\\
=&-\frac{1}{8}(\Lambda_{n}^{-1}a_{n}^{-2}b_{n}^{-2})^{2}(1+\beta_{n}i)\alpha_{2n-1}^{-1}-\frac{1}{16i}\Lambda_{n}^{-1}a_{n}^{-2}b_{n}^{-2}\alpha_{2n-1}^{-1}[a_{n}^{-2}(1+\beta_{n}^{2})+b_{n}^{-2}]\nonumber\\
&+\frac{1}{16i}\Lambda_{n}^{-1}a_{n}^{-2}b_{n}^{-2}[b_{n}^{-2}-a_{n}^{-2}(1+\beta_{n}^{2})]\nonumber\\
=&\frac{1}{16i}\Lambda_{n}^{-1}a_{n}^{-2}b_{n}^{-2}\alpha_{2n-1}^{-1}\Big\{\kappa_{2n}^{-2}a_{n}^{-2}b_{n}^{-2}(1+\beta_{n}i)-[a_{n}^{-2}(1+\beta_{n}^{2})+b_{n}^{-2}]\nonumber\\
&+\alpha_{2n-1}[b_{n}^{-2}-a_{n}^{-2}(1+\beta_{n}^{2})]\Big\}\nonumber\\
=&-\frac{1}{8i}\Lambda_{n}^{-1}a_{n}^{-4}b_{n}^{-2}\alpha_{2n-1}^{-1}\beta_{n}(\beta_{n}-i).\nonumber
\end{align}
Thus the proof is complete.
\end{proof}

\begin{cor}
Let $B, C, D, E, F$ be the coefficient matrices as above, then $B, C, D, E$ are invertible for $\beta_{n}\in \mathbb{R}$, whereas $F$ is invertible for $\beta_{n}\in \mathbb{R}\setminus\{0\}$.
\end{cor}

By Theorems 5.2 and 5.3 as well as Corollary 5.6, we have
\begin{thm}
Let $A, B, C, D, E, F$ be the coefficient matrices as above, then
\begin{align}
&A^{-1}\left(
                \begin{array}{c}
                  \imath_{n+1}-\frac{1}{4}a_{n}^{-2}(1-\beta_{n}i)\overline{\alpha}_{2n}^{-1}\Big(\kappa_{2n}^{-2}-\kappa_{2n+1}^{-2}\Big) \vspace{1mm}\\
                  \jmath_{n+1}-\frac{1}{4}b_{n}^{-2}\overline{\alpha}_{2n}^{-1}\Big(\kappa_{2n}^{-2}-\kappa_{2n+1}^{-2}\Big)i \\
                \end{array}
              \right)\nonumber\\
              =&B^{-1}\left(
                \begin{array}{c}
                  \varsigma_{n+1}+\frac{1}{4i}a_{n}^{-2}(1-\beta_{n}i)\overline{\alpha}_{2n}^{-1}\Big(\kappa_{2n}^{-2}-\kappa_{2n+1}^{-2}\Big) \vspace{1mm}\\
                  \zeta_{n+1}+\frac{1}{4}b_{n}^{-2}\overline{\alpha}_{2n}^{-1}\Big(\kappa_{2n}^{-2}-\kappa_{2n+1}^{-2}\Big) \\
                \end{array}
              \right)\nonumber\\
              =&C^{-1}\left(
                \begin{array}{c}
                  \imath_{n+1}-\frac{1}{4}a_{n}^{-2}(1-\beta_{n}i)\overline{\alpha}_{2n}^{-1}\Big(\kappa_{2n}^{-2}-\kappa_{2n+1}^{-2}\Big) \vspace{1mm}\\
                  \varsigma_{n+1}+\frac{1}{4i}a_{n}^{-2}(1-\beta_{n}i)\overline{\alpha}_{2n}^{-1}\Big(\kappa_{2n}^{-2}-\kappa_{2n+1}^{-2}\Big) \\
                \end{array}
              \right)\nonumber\\
              =&D^{-1}\left(
                \begin{array}{c}
                  \jmath_{n+1}-\frac{1}{4}b_{n}^{-2}\overline{\alpha}_{2n}^{-1}\Big(\kappa_{2n}^{-2}-\kappa_{2n+1}^{-2}\Big)i \vspace{1mm}\\
                  \zeta_{n+1}+\frac{1}{4}b_{n}^{-2}\overline{\alpha}_{2n}^{-1}\Big(\kappa_{2n}^{-2}-\kappa_{2n+1}^{-2}\Big) \\
                \end{array}
              \right)\nonumber\\
              =&E^{-1}\left(
                \begin{array}{c}
                  \imath_{n+1}-\frac{1}{4}a_{n}^{-2}(1-\beta_{n}i)\overline{\alpha}_{2n}^{-1}\Big(\kappa_{2n}^{-2}-\kappa_{2n+1}^{-2}\Big) \vspace{1mm}\\
                  \zeta_{n+1}+\frac{1}{4}b_{n}^{-2}\overline{\alpha}_{2n}^{-1}\Big(\kappa_{2n}^{-2}-\kappa_{2n+1}^{-2}\Big) \\
                \end{array}
              \right)
\end{align}
for $\beta_{n}\in \mathbb{R}$. Moreover, any term in the above identities is equal to
\begin{equation}
F^{-1}\left(
                \begin{array}{c}
                  \jmath_{n+1}-\frac{1}{4}b_{n}^{-2}\overline{\alpha}_{2n}^{-1}\Big(\kappa_{2n}^{-2}-\kappa_{2n+1}^{-2}\Big)i \vspace{1mm}\\
                  \varsigma_{n+1}+\frac{1}{4i}a_{n}^{-2}(1-\beta_{n}i)\overline{\alpha}_{2n}^{-1}\Big(\kappa_{2n}^{-2}-\kappa_{2n+1}^{-2}\Big) \\
                \end{array}
              \right)
\end{equation}
as $\beta_{n}\neq 0$.
\end{thm}

Let
\begin{equation}
\gamma=\left(
                \begin{array}{c}
                  \frac{1}{4i}a_{n}^{-2}(1-\beta_{n}i)\overline{\alpha}_{2n}^{-1}\Big(\kappa_{2n}^{-2}-\kappa_{2n+1}^{-2}\Big) \vspace{1mm}\\
                  \frac{1}{4}b_{n}^{-2}\overline{\alpha}_{2n}^{-1}\Big(\kappa_{2n}^{-2}-\kappa_{2n+1}^{-2}\Big) \\
                \end{array}
              \right)
\end{equation}
and
\begin{equation}
\eta=\left(
                \begin{array}{c}
                  \frac{1}{4}a_{n}^{-2}(1-\beta_{n}i)\overline{\alpha}_{2n}^{-1}\Big(\kappa_{2n}^{-2}-\kappa_{2n+1}^{-2}\Big) \vspace{1mm}\\
                  \frac{1}{4}b_{n}^{-2}\overline{\alpha}_{2n}^{-1}\Big(\kappa_{2n}^{-2}-\kappa_{2n+1}^{-2}\Big)i \\
                \end{array}
              \right),
\end{equation}
then by the first identity in (5.25),
\begin{equation}
B\left(
   \begin{array}{c}
     \imath_{n+1} \\
     \jmath_{n+1} \\
   \end{array}
 \right)-A\left(
   \begin{array}{c}
     \varsigma_{n+1} \\
     \zeta_{n+1} \\
   \end{array}
 \right)=A\gamma+B\eta,
\end{equation}
where $A$ and $B$ are given by (5.16) and (5.19) respectively.

By (2.23), we have

\begin{align}
\left(
  \begin{array}{c}
    \alpha_{2n} \vspace{1mm}\\
    \overline{\alpha}_{2n} \\
  \end{array}
\right)=P\left(
          \begin{array}{c}
            \imath_{n+1} \vspace{1mm}\\
            \jmath_{n+1} \\
          \end{array}
        \right)+Q
        \left(
          \begin{array}{c}
            \varsigma_{n+1} \vspace{1mm}\\
            \zeta_{n+1} \\
          \end{array}
        \right),
\end{align}
where
\begin{equation}
P=\left(
          \begin{array}{cc}
            \frac{1+\beta_{n}i}{2} & -\frac{i}{2} \vspace{1mm}\\
            \frac{1-\beta_{n}i}{2} & \frac{i}{2} \\
          \end{array}
        \right)\,\,\,\mbox{and}\,\,\,Q=\left(
          \begin{array}{cc}
            \frac{\beta_{n}-i}{2} & -\frac{1}{2} \vspace{1mm}\\
            \frac{\beta_{n}+i}{2} & -\frac{1}{2} \\
          \end{array}
        \right).
\end{equation}
Denote
\begin{equation}
D=\left(
                \begin{array}{cc}
                  P & Q \\
                  B & -A \\
                \end{array}
              \right)\,\,\,\mbox{and}\,\,\,\alpha=\left(
         \begin{array}{c}
           \alpha_{2n} \\
           \overline{\alpha}_{2n} \\
         \end{array}
       \right),
\end{equation}
then we have the following theorem.
\begin{thm}
\begin{equation}
\left(
         \begin{array}{c}
           \imath_{n+1} \\
           \jmath_{n+1} \\
           \varsigma_{n+1} \\
           \zeta_{n+1} \\
         \end{array}
       \right)=D^{-1}\left(
                 \begin{array}{c}
                   \alpha \\
                   A\gamma+B\eta \\
                 \end{array}
               \right),
\end{equation}
where $D^{-1}$ is the inverse matrix of $D$.
\end{thm}

\begin{proof}
By (5.29) and (5.30),
\begin{equation}
D\left(
         \begin{array}{c}
           \imath_{n+1} \\
           \jmath_{n+1} \\
           \varsigma_{n+1} \\
           \zeta_{n+1} \\
         \end{array}
       \right)=\left(
  \begin{array}{cc}
    P & Q \\
    B & -A \\
  \end{array}
\right)\left(
         \begin{array}{c}
           \imath_{n+1} \\
           \jmath_{n+1} \\
           \varsigma_{n+1} \\
           \zeta_{n+1} \\
         \end{array}
       \right)=\left(
                 \begin{array}{c}
                   \alpha \\
                   A\gamma+B\eta \\
                 \end{array}
               \right).
\end{equation}
So it is enough to show that $|D|\neq0$ in order to get (5.33).  However, it follows from Laplace theorem by simple calculations on some of its subdeterminants. More precisely,
\begin{itemize}
  \item [Case 1.] \begin{equation}
  |P|=\left|
          \begin{array}{cc}
            \frac{1+\beta_{n}i}{2} & -\frac{i}{2} \vspace{1mm}\\
            \frac{1-\beta_{n}i}{2} & \frac{i}{2} \\
          \end{array}
        \right|=\frac{i}{2}\frac{1+\beta_{n}i}{2}+\frac{i}{2}\frac{1-\beta_{n}i}{2}=\frac{i}{2}
  \end{equation}
  and
  \begin{equation}
  (-1)^{1+2+1+2}|-A|=|A|.
  \end{equation}
  \item [Case 2.] \begin{equation}
  |Q|=\left|
          \begin{array}{cc}
            \frac{\beta_{n}-i}{2} & -\frac{1}{2} \vspace{1mm}\\
            \frac{\beta_{n}+i}{2} & -\frac{1}{2} \\
          \end{array}
        \right|=-\frac{1}{2}\frac{\beta_{n}-i}{2}+\frac{1}{2}\frac{\beta_{n}+i}{2}=\frac{i}{2}
  \end{equation}
  and
  \begin{equation}
  (-1)^{1+2+3+4}|B|=|B|.
  \end{equation}
  \item [Case 3.] \begin{equation}
  \left|
          \begin{array}{cc}
            \frac{1+\beta_{n}i}{2} & \frac{\beta_{n}-i}{2} \vspace{1mm}\\
            \frac{1-\beta_{n}i}{2} & \frac{\beta_{n}+i}{2} \\
          \end{array}
        \right|=\frac{1+\beta_{n}i}{2}\frac{\beta_{n}+i}{2}-\frac{\beta_{n}-i}{2}\frac{1-\beta_{n}i}{2}=\frac{1+\beta_{n}^{2}}{2}i
        \end{equation}
        and
        \begin{align}
  &(-1)^{1+2+1+3}\left|
          \begin{array}{cc}
            -\frac{1}{4i}\Big[\Lambda_{n}^{-1}a_{n}^{-2}b_{n}^{-2}\alpha_{2n-1}^{-1}i+a_{n}^{-2}(1-\beta_{n}i)\Big] & \frac{1}{4}\Big[\Lambda_{n}^{-1}a_{n}^{-2}b_{n}^{-2}\alpha_{2n-1}^{-1}i+a_{n}^{-2}(1-\beta_{n}i)\Big] \vspace{1mm}\\
            -\frac{1}{4i}\Big[\Lambda_{n}^{-1}a_{n}^{-2}b_{n}^{-2}(1+\beta_{n}i)\alpha_{2n-1}^{-1}+b_{n}^{-2}i\Big] & \frac{1}{4}\Big[\Lambda_{n}^{-1}a_{n}^{-2}b_{n}^{-2}(1+\beta_{n}i)\alpha_{2n-1}^{-1}+b_{n}^{-2}i\Big] \\
          \end{array}
        \right|\nonumber\\
        &=0.
        \end{align}

  \item [Case 4.] \begin{equation}\left|
          \begin{array}{cc}
            \frac{1+\beta_{n}i}{2} & -\frac{1}{2} \vspace{1mm}\\
            \frac{1-\beta_{n}i}{2} & -\frac{1}{2} \\
          \end{array}
        \right|=-\frac{1}{2}\frac{1+\beta_{n}i}{2}+\frac{1}{2}\frac{1-\beta_{n}i}{2}=-\frac{\beta_{n}}{2}i
  \end{equation}
  and
  \begin{align}
  &(-1)^{1+2+1+4}\left|
          \begin{array}{cc}
            -\frac{1}{4i}\Big[\Lambda_{n}^{-1}a_{n}^{-2}b_{n}^{-2}\alpha_{2n-1}^{-1}i+a_{n}^{-2}(1-\beta_{n}i)\Big] & -\frac{1}{4}\Lambda_{n}^{-1}a_{n}^{-2}b_{n}^{-2}\Big(\alpha_{2n-1}^{-1}+1\Big)i \vspace{1mm}\\
            -\frac{1}{4i}\Big[\Lambda_{n}^{-1}a_{n}^{-2}b_{n}^{-2}(1+\beta_{n}i)\alpha_{2n-1}^{-1}+b_{n}^{-2}i\Big] & -\frac{1}{4}\Lambda_{n}^{-1}a_{n}^{-2}b_{n}^{-2}(1+\beta_{n}i)\Big(\alpha_{2n-1}^{-1}+1\Big) \\
          \end{array}
        \right|\nonumber\\
        =&-|A|i.
  \end{align}
  \item [Case 5.] \begin{equation}
  \left|
          \begin{array}{cc}
            -\frac{i}{2} & \frac{\beta_{n}-i}{2} \vspace{1mm}\\
            \frac{i}{2} & \frac{\beta_{n}+i}{2} \\
          \end{array}
        \right|=-\frac{i}{2}\frac{\beta_{n}+i}{2}-\frac{\beta_{n}-i}{2}\frac{i}{2}=-\frac{\beta_{n}}{2}i
        \end{equation}
        and
        \begin{align}
  &(-1)^{1+2+2+3}\left|
          \begin{array}{cc}
            \frac{1}{4}\Lambda_{n}^{-1}a_{n}^{-2}b_{n}^{-2}\Big(\alpha_{2n-1}^{-1}-1\Big) & \frac{1}{4}\Big[\Lambda_{n}^{-1}a_{n}^{-2}b_{n}^{-2}\alpha_{2n-1}^{-1}i+a_{n}^{-2}(1-\beta_{n}i)\Big] \vspace{1mm}\\
            \frac{1}{4i}\Lambda_{n}^{-1}a_{n}^{-2}b_{n}^{-2}(1+\beta_{n}i)\Big(\alpha_{2n-1}^{-1}-1\Big) & \frac{1}{4}\Big[\Lambda_{n}^{-1}a_{n}^{-2}b_{n}^{-2}(1+\beta_{n}i)\alpha_{2n-1}^{-1}+b_{n}^{-2}i\Big] \\
          \end{array}
        \right|\nonumber\\
        =&-|B|i.
        \end{align}
  \item [Case 6.]
\end{itemize} \begin{equation}
  \left|
          \begin{array}{cc}
            -\frac{i}{2} & -\frac{1}{2} \vspace{1mm}\\
            \frac{i}{2} & -\frac{1}{2} \\
          \end{array}
        \right|=\frac{i}{2}\frac{1}{2}+\frac{i}{2}\frac{1}{2}=\frac{i}{2}
        \end{equation}
        and
        \begin{align}
  &(-1)^{1+2+2+4}\left|
          \begin{array}{cc}
            \frac{1}{4}\Lambda_{n}^{-1}a_{n}^{-2}b_{n}^{-2}\Big(\alpha_{2n-1}^{-1}-1\Big) & -\frac{1}{4}\Lambda_{n}^{-1}a_{n}^{-2}b_{n}^{-2}\Big(\alpha_{2n-1}^{-1}+1\Big)i \vspace{1mm}\\
            \frac{1}{4i}\Lambda_{n}^{-1}a_{n}^{-2}b_{n}^{-2}(1+\beta_{n}i)\Big(\alpha_{2n-1}^{-1}-1\Big) & -\frac{1}{4}\Lambda_{n}^{-1}a_{n}^{-2}b_{n}^{-2}(1+\beta_{n}i)\Big(\alpha_{2n-1}^{-1}+1\Big) \\
          \end{array}
        \right|=0.
        \end{align}
So
  \begin{align}
  |D|=&\frac{i}{2}(|A|+|B|)-\frac{\beta_{n}}{2}i(-|A|i-|B|i)=(|A|+|B|)\frac{i-\beta_{n}}{2}\nonumber\\
  =&(|A|+|B|)i\frac{1+\beta_{n}i}{2}=-\frac{1}{8}a_{n}^{-2}b_{n}^{-2}(1+\beta_{n}i)\neq 0.
  \end{align}
Thus
\begin{equation}
\left(
         \begin{array}{c}
           \imath_{n+1} \\
           \jmath_{n+1} \\
           \varsigma_{n+1} \\
           \zeta_{n+1} \\
         \end{array}
       \right)=D^{-1}\left(
                 \begin{array}{c}
                   \alpha \\
                   A\gamma+B\eta \\
                 \end{array}
               \right).\nonumber\vspace{-5mm}
\end{equation}
\end{proof}

\begin{rem}
Similar as in Remark 2.10, we can directly obtain (5.33) by using (2.21)-(2.24) and (2.37) together.
\end{rem}

With the above preliminaries, we get the following strong Favard theorem.

\begin{thm} Let
$\{(a_{n}^{(0)}, b_{n}^{(0)},\beta_{n}^{(0)},\imath_{n}^{(0)},\jmath_{n}^{(0)},\varsigma_{n}^{(0)},\zeta_{n}^{(0)})\}$ with
$a_{0}^{(0)},b_{0}^{(0)}=1$ and $\beta_{0}^{(0)}=0$ be a system of
seven-tuples of real numbers satisfying
\begin{itemize}
  \item [(1)] \begin{align}
&0<(\imath_{n+1}^{(0)}+\beta_{n}^{(0)}\varsigma_{n+1}^{(0)}-\zeta_{n+1}^{(0)})^{2}+(\jmath_{n+1}^{(0)}-\imath_{n+1}^{(0)}\beta_{n}^{(0)}
+\varsigma_{n+1}^{(0)})^{2}\nonumber\\
=&1-\frac{1}{4}[(a_{n}^{(0)})^{2}+(b_{n}^{(0)})^{2}(1+(\beta_{n}^{(0)})^{2})]
[(a_{n+1}^{(0)})^{2}+(b_{n+1}^{(0)})^{2}(1+(\beta_{n+1}^{(0)})^{2})]<1;
\end{align}
  \item [(2)]
  \begin{equation}
  \beta_{n}^{(0)}\neq0\,\,\,\,\mbox{or}\,\,\,\,\frac{(a_{n}^{(0)})^{2}}{(b_{n}^{(0)})^{2}}+(\beta_{n}^{(0)})^{2}\neq1
  \end{equation}
\end{itemize}
for $n\in \mathbb{N}$ with $a_{n}^{(0)},b_{n}^{(0)}>0$ for $n\in \mathbb{N}\cup\{0\}$,
then there exists a unique nontrivial probability measure $d\mu$ on
$\partial \mathbb{D}$ such that $a_{n}(d\mu)=a_{n}^{(0)}$,
$b_{n}(d\mu)=b_{n}^{(0)}$, $\beta_{n}(d\mu)=\beta_{n}^{(0)}$, $\imath_{n}(d\mu)=\imath_{n}^{(0)}$, $\jmath_{n}(d\mu)=\jmath_{n}^{(0)}$, $\varsigma_{n}(d\mu)=\varsigma_{n}^{(0)}$ and $\zeta_{n}(d\mu)=\zeta_{n}^{(0)}$,
where $a_{n}(d\mu),b_{n}(d\mu),\beta_{n}(d\mu),\imath_{n}(d\mu),\jmath_{n}(d\mu),\varsigma_{n}(d\mu),\zeta_{n}(d\mu)$ are associated
coefficients of $d\mu$ defined by (2.11)-(2.15).
\end{thm}

\begin{proof}
Let \begin{equation}
\kappa_{2n}^{(0)}=\frac{1}{2}\Big[(a_{n}^{(0)})^{-2}\big(1+(\beta_{n}^{(0)})^{2}\big)+(b_{n}^{(0)})^{-2}\Big]^{\frac{1}{2}},
\end{equation}
\begin{equation}
\kappa_{2n+1}^{(0)}=\Big[(a_{n+1}^{(0)})^{2}+(b_{n+1}^{(0)})^{2}\big(1+(\beta_{n+1}^{(0)})^{2}\big)\Big]^{-\frac{1}{2}},
\end{equation}
\begin{equation}
\alpha_{2n-1}^{(0)}=\frac{1}{4}(\kappa_{2n}^{(0)})^{-2}\Big[(b_{n}^{(0)})^{-2}-(a_{n}^{(0)})^{-2}\big(1-(\beta_{n}^{(0)})^{2}\big)\Big]
-\frac{1}{2}(\kappa_{2n}^{(0)})^{-2}(a_{n}^{(0)})^{-2}
(\beta_{n}^{(0)})i
\end{equation}
and
\begin{equation}
\alpha_{2n}^{(0)}=\frac{1}{2}(\imath_{n+1}^{(0)}+\beta_{n}^{(0)}\varsigma_{n+1}^{(0)}-\zeta_{n+1}^{(0)})-\frac{i}{2}(\jmath_{n+1}^{(0)}-\imath_{n+1}^{(0)}\beta_{n}^{(0)}
+\varsigma_{n+1}^{(0)})
\end{equation}
for $n\in \mathbb{N}\cup\{0\}$, then by (5.48), Verblunsky theorem and Theorem 3.1, there exists a unique nontrivial probability measure $d\mu$ on $\partial \mathbb{D}$ such that
\begin{equation}\alpha_{n}(d\mu)=\alpha_{n}^{(0)}\,\,\,\mbox{and}\,\,\,\kappa_{n}(d\mu)=\kappa_{n}^{(0)}\end{equation}
as well as
\begin{equation}a_{n}(d\mu)=a_{n}^{(0)},\,\,b_{n}(d\mu)=b_{n}^{(0)}\,\,\,\mbox{and}\,\,\,\beta_{n}(d\mu)=\beta_{n}^{(0)}\end{equation} for $n\in \mathbb{N}\cup\{0\}$.

On one hand, noting (5.49), by Theorem 5.8, we have
\begin{equation}
\left(
         \begin{array}{c}
           \imath_{n+1}(d\mu) \vspace{1mm}\\
           \jmath_{n+1}(d\mu) \vspace{1mm}\\
           \varsigma_{n+1}(d\mu) \vspace{1mm}\\
           \zeta_{n+1}(d\mu) \\
         \end{array}
       \right)=D(d\mu)^{-1}\left(
                 \begin{array}{c}
                   \alpha(d\mu) \\
                   A(d\mu)\gamma(d\mu)+B(d\mu)\eta(d\mu) \\
                 \end{array}
               \right)
\end{equation}
for $n\in \mathbb{N}\cup\{0\}$, where $A(d\mu), B(d\mu), \gamma(d\mu), \eta(d\mu), D(d\mu), \alpha(d\mu)$ are respectively given by (5.16), (5.19), (5.27), (5.28) and (5.32) with $a_{n}(d\mu), b_{n}(d\mu),\beta_{n}(d\mu), \alpha_{n}(d\mu),$ $\kappa_{n}(d\mu), \Lambda_{n}(d\mu)$ replacing $a_{n}, b_{n},\beta_{n}, \alpha_{n}, \kappa_{n}, \Lambda_{n}$.

On the other hand, as in Remark 5.9, by directly invoking (5.50)-(5.53) and
\begin{equation}\Lambda_{n}^{(0)}=-\frac{1}{2}\Big[(a_{n}^{(0)})^{-2}\big(1+(\beta_{n}^{(0)})^{2}\big)+(b_{n}^{(0)})^{-2}\Big]i,\end{equation}
we have
\begin{equation}
\left(
         \begin{array}{c}
           \imath_{n+1}^{(0)} \vspace{1mm}\\
           \jmath_{n+1}^{(0)} \vspace{1mm}\\
           \varsigma_{n+1}^{(0)} \vspace{1mm}\\
           \zeta_{n+1}^{(0)} \vspace{1mm}\\
         \end{array}
       \right)=(D^{(0)})^{-1}\left(
                 \begin{array}{c}
                   \alpha^{(0)} \\
                   A^{(0)}\gamma^{(0)}+B^{(0)}\eta^{(0)} \\
                 \end{array}
               \right)
\end{equation}
for $n\in \mathbb{N}\cup\{0\}$, where $A^{(0)}, B^{(0)}, \gamma^{(0)}, \eta^{(0)}, D^{(0)}, \alpha^{(0)}$ are respectively given by (5.16), (5.19), (5.27), (5.28) and (5.32) with $a_{n}^{(0)}, b_{n}^{(0)},\beta_{n}^{(0)}, \alpha_{n}^{(0)},$ $\kappa_{n}^{(0)}, \Lambda_{n}^{(0)}$ replacing $a_{n}, b_{n},\beta_{n}, \alpha_{n}, \kappa_{n}, \Lambda_{n}$.

In terms of (2.37), (5.16), (5.19), (5.27), (5.28), (5.32), (5.54), (5.55) and (5.57), one can easily find that
\begin{equation*}D(d\mu)=D^{(0)},A(d\mu)=A^{(0)}, B(d\mu)=B^{(0)}, \alpha(d\mu)=\alpha^{(0)},\gamma(d\mu)=\gamma^{(0)},\eta(d\mu)=\eta^{(0)}.\end{equation*}
Thus, by (5.56) and (5.58),
\begin{equation*}\imath_{n}(d\mu)=\imath_{n}^{(0)}, \jmath_{n}(d\mu)=\jmath_{n}^{(0)}, \varsigma_{n}(d\mu)=\varsigma_{n}^{(0)}\,\,\,\mbox{and}\,\,\, \zeta_{n}(d\mu)=\zeta_{n}^{(0)}\end{equation*}
for $n\in \mathbb{N}$.
\end{proof}

Except for the above strong Favard theorem, by Theorem 5.5 and (5.17), we also have the following result on the determinants of the coefficient matrices $A, B, C, D, E$ and $F$.

\begin{thm}
Let $A, B, C, D, E, F$ be the coefficient matrices as above, then
\begin{equation}
|A||B|-|C||D|(1+\beta_{n}i)+|E||F|=0.
\end{equation}
\end{thm}

\begin{proof}
Note that
\begin{align}
FE^{-1}=&|E|^{-1}S,
\end{align}
where
\begin{align}
S=&\left(\begin{array}{cc}
      \frac{1}{4}\Lambda_{n}^{-1}a_{n}^{-2}b_{n}^{-2}(1+\beta_{n}i)\Big(\alpha_{2n-1}^{-1}+1\Big) & -\frac{1}{4}[\Lambda_{n}^{-1}a_{n}^{-2}b_{n}^{-2}(1+\beta_{n}i)\alpha_{2n-1}^{-1}+b_{n}^{-2}i]\vspace{2mm}\\
      \frac{1}{4}\Lambda_{n}^{-1}a_{n}^{-2}b_{n}^{-2}\Big(\alpha_{2n-1}^{-1}-1\Big) & -\frac{1}{4i}[\Lambda_{n}^{-1}a_{n}^{-2}b_{n}^{-2}\alpha_{2n-1}^{-1}i+a_{n}^{-2}(1-\beta_{n}i)] \\
    \end{array}
\right)\nonumber\\
&\left(\begin{array}{cc}
      -\frac{1}{4i}[\Lambda_{n}^{-1}a_{n}^{-2}b_{n}^{-2}(1+\beta_{n}i)\alpha_{2n-1}^{-1}+b_{n}^{-2}i] & \frac{1}{4}[\Lambda_{n}^{-1}a_{n}^{-2}b_{n}^{-2}\alpha_{2n-1}^{-1}i+a_{n}^{-2}(1-\beta_{n}i)]\vspace{2mm}\\
      -\frac{1}{4i}\Lambda_{n}^{-1}a_{n}^{-2}b_{n}^{-2}(1+\beta_{n}i)\Big(\alpha_{2n-1}^{-1}-1\Big) & \frac{1}{4}\Lambda_{n}^{-1}a_{n}^{-2}b_{n}^{-2}\Big(\alpha_{2n-1}^{-1}+1\Big)i\\
    \end{array}
\right)\nonumber\\
\triangleq&\left(
    \begin{array}{cc}
      s_{11} & s_{12} \\
      s_{21} & s_{22} \\
    \end{array}
  \right)
\end{align}
with
\begin{align}
s_{11}=&\frac{1}{4}\Lambda_{n}^{-1}a_{n}^{-2}b_{n}^{-2}(1+\beta_{n}i)\Big(\alpha_{2n-1}^{-1}+1\Big)
\Big\{-\frac{1}{4i}[\Lambda_{n}^{-1}a_{n}^{-2}b_{n}^{-2}(1+\beta_{n}i)\alpha_{2n-1}^{-1}+b_{n}^{-2}i]\Big\}\nonumber\\
&-\frac{1}{4}[\Lambda_{n}^{-1}a_{n}^{-2}b_{n}^{-2}(1+\beta_{n}i)\alpha_{2n-1}^{-1}+b_{n}^{-2}i]\Big\{-\frac{1}{4i}\Lambda_{n}^{-1}a_{n}^{-2}b_{n}^{-2}(1+\beta_{n}i)\Big(\alpha_{2n-1}^{-1}-1\Big) \Big\}\nonumber\\
=&-\frac{1}{8i}\Lambda_{n}^{-1}a_{n}^{-2}b_{n}^{-2}(1+\beta_{n}i)[\Lambda_{n}^{-1}a_{n}^{-2}b_{n}^{-2}(1+\beta_{n}i)\alpha_{2n-1}^{-1}+b_{n}^{-2}i]\nonumber\\
=&-\frac{1}{8}\Lambda_{n}^{-1}a_{n}^{-2}b_{n}^{-4}\alpha_{2n-1}^{-1}(1+\beta_{n}i)=|D|(1+\beta_{n}i),
\end{align}

\begin{align}
s_{12}=&\frac{1}{4}\Lambda_{n}^{-1}a_{n}^{-2}b_{n}^{-2}(1+\beta_{n}i)\Big(\alpha_{2n-1}^{-1}+1\Big)
\Big\{\frac{1}{4}[\Lambda_{n}^{-1}a_{n}^{-2}b_{n}^{-2}\alpha_{2n-1}^{-1}i+a_{n}^{-2}(1-\beta_{n}i)]\Big\}\nonumber\\
&-\frac{1}{4}[\Lambda_{n}^{-1}a_{n}^{-2}b_{n}^{-2}(1+\beta_{n}i)\alpha_{2n-1}^{-1}+b_{n}^{-2}i]\Big\{\frac{1}{4}\Lambda_{n}^{-1}a_{n}^{-2}b_{n}^{-2}\Big(\alpha_{2n-1}^{-1}+1\Big)i\Big\}\nonumber\\
=&\frac{1}{16}\Lambda_{n}^{-1}a_{n}^{-2}b_{n}^{-2}[a_{n}^{-2}(1+\beta_{n}^{2})+b_{n}^{-2}]\Big(\alpha_{2n-1}^{-1}+1\Big)\nonumber\\
=&\frac{1}{4}\Lambda_{n}^{-1}a_{n}^{-2}b_{n}^{-2}\kappa_{2n}^{2}\Big(\alpha_{2n-1}^{-1}+1\Big)=|A|,
\end{align}

\begin{align}
s_{21}=&\frac{1}{4}\Lambda_{n}^{-1}a_{n}^{-2}b_{n}^{-2}\Big(\alpha_{2n-1}^{-1}-1\Big)
\Big\{-\frac{1}{4i}[\Lambda_{n}^{-1}a_{n}^{-2}b_{n}^{-2}(1+\beta_{n}i)\alpha_{2n-1}^{-1}+b_{n}^{-2}i]\Big\}\nonumber\\
&-\frac{1}{4i}[\Lambda_{n}^{-1}a_{n}^{-2}b_{n}^{-2}\alpha_{2n-1}^{-1}i+a_{n}^{-2}(1-\beta_{n}i)]\Big\{-\frac{1}{4i}\Lambda_{n}^{-1}a_{n}^{-2}b_{n}^{-2}(1+\beta_{n}i)\Big(\alpha_{2n-1}^{-1}-1\Big)\Big\}\nonumber\\
=&-\frac{1}{16}\Lambda_{n}^{-1}a_{n}^{-2}b_{n}^{-2}[a_{n}^{-2}(1+\beta_{n}^{2})+b_{n}^{-2}]\Big(\alpha_{2n-1}^{-1}-1\Big)\nonumber\\
=&-\frac{1}{4}\Lambda_{n}^{-1}a_{n}^{-2}b_{n}^{-2}\kappa_{2n}^{2}\Big(\alpha_{2n-1}^{-1}-1\Big)=|B|,
\end{align}
and
\begin{align}
s_{22}=&\frac{1}{4}\Lambda_{n}^{-1}a_{n}^{-2}b_{n}^{-2}\Big(\alpha_{2n-1}^{-1}-1\Big)
\Big\{\frac{1}{4}[\Lambda_{n}^{-1}a_{n}^{-2}b_{n}^{-2}\alpha_{2n-1}^{-1}i+a_{n}^{-2}(1-\beta_{n}i)]\Big\}\nonumber\\
&-\frac{1}{4i}[\Lambda_{n}^{-1}a_{n}^{-2}b_{n}^{-2}\alpha_{2n-1}^{-1}i+a_{n}^{-2}(1-\beta_{n}i)]\Big\{\frac{1}{4}\Lambda_{n}^{-1}a_{n}^{-2}b_{n}^{-2}\Big(\alpha_{2n-1}^{-1}+1\Big)i\Big\}\nonumber\\
=&-\frac{1}{8}\Lambda_{n}^{-1}a_{n}^{-2}b_{n}^{-2}[\Lambda_{n}^{-1}a_{n}^{-2}b_{n}^{-2}\alpha_{2n-1}^{-1}i+a_{n}^{-2}(1-\beta_{n}i)]\nonumber\\
=&\frac{1}{8}\Lambda_{n}^{-1}a_{n}^{-4}b_{n}^{-2}\alpha_{2n-1}^{-1}(1+\beta_{n}i)=|C|.
\end{align}
The last steps in the above identities for all the entries of $S$ hold by (5.17) and Theorem 5.5.
So
\begin{align}
FE^{-1}=|E|^{-1}\left(
                  \begin{array}{cc}
                    |D|(1+\beta_{n}i) & |A| \\
                    |B| & |C| \\
                  \end{array}
                \right).
\end{align}
Thus (5.59) follows immediately.
\end{proof}

\begin{rem}
In fact, one can get (5.59) by directly invoking (5.20)-(5.24) and (5.17). However, by this means, (5.66) could not be seen.
\end{rem}

\bibliographystyle{amsplain}

\begin{thebibliography}{10}


\bibitem{bdj}J. Baik, P. Deift, and K. Johansson, {\it On the distribution of the length of the longest increasing subsequence
of random permutations}, J. Am. Math. Soc. \textbf{12} (1999), 1119-1179.

\bibitem{deift}P. Deift, {\it Orthogonal Polynomials and Random Matrices: A Riemann-Hilbert Approach}, Courant Lecture Notes in Math.,
Vol. 3, Amer. Math. Soc., Providence, RI, 2000.

\bibitem{deift1}P. Deift, {\it Riemann-Hilbert Problems},  Random matrices, 1-40, IAS/Park City Math. Ser., 26, Amer. Math. Soc., Providence, RI, 2019.

\bibitem{dz}P. Deift, X. Zhou, {\it A steepest descent method for oscillatory Riemann-Hilbert problems. Asymptotics
for the MKdV Equation}, Ann. Math. \textbf{137} (1993), 295-368.

\bibitem{d}Z. Du, {\it A semi-conjugate matrix boundary value problem for general orthogonal
polynomials on an arbitrary smooth Jordan curve}, Acta Math. Sci. \textbf{28B} (2008), 401-407.

\bibitem{dd06} Z. Du and J. Du, Riemann-Hilbert approach to strong asymptotics for orthogonal polynomials on the unit circle, Chinese Ann. Math. Ser A \textbf{27} (2006),
701-718; also in Chinese J. Contemp. Math. \textbf{27} (2006), 443-462.

\bibitem{dd08}Z. Du and J. Du, {\it Orthogonal trigonometric
polynomials: Riemann-Hilbert analysis and Relations with OPUC},
Asymptotic Analysis \textbf{79} (2012) 87-132.

\bibitem{enzg}T. Erd\'elyi, P. Nevai, J. Zhang and J. Geronimo, {\it A simple proof of ``Favard's theorem" on the unit circle},
Atti Sem. Mat. Fis. Univ. Modena \textbf{39} (1991), 551-556.

\bibitem{fik} A. Fokas, A. Its, and A. Kitaev, {\it The isomonodromy approach to matrix models in 2D quantum
gravity}, Commun. Math. Phys. \textbf{147} (1992), 395-430.

\bibitem{ma}F. Marcell\'an and R. \'Alvarez-Nodarse, {\it On the ``Favard theorem" and its extensions},
J. Comp. Appl. Math. \textbf{127} (2001), 231-254.

\bibitem{mehta}M. Mehta, {\it Random Matrices}, 2nd ed.,
Academic, Boston, 1991.

\bibitem{sim1}B. Simon, {\it Orthogonal Polynomials on the Unit Circle, Part 1: Classical Theory},
AMS Colloquium Series, Vol. 54, American Mathematical Society,
Providence RI, 2005.

\bibitem{sim2}B. Simon, {\it Orthogonal Polynomials on the Unit Circle, Part 2: Spectral Theory},
AMS Colloquium Series, Vol. 54, American Mathematical Society,
Providence RI, 2005.

\bibitem{sze}G. Szeg\"o, {\it Orthogonal Polynomials}, Amer. Math. Soc.
Colloq. Publ., Vol. 23, 4th ed., American Mathematical Society, Providence
R.I., 1975.

\end{thebibliography}

\end{document}